\documentclass{amsart}
\usepackage{amsfonts}
\newtheorem{lem}{Lemma}

\DeclareMathOperator{\Tr}{Tr}

\DeclareMathOperator{\diag}{diag}
\newcommand{\abs}[1]{\left\vert#1\right\vert}

\newtheorem{thm}{Theorem}

\let\f=\frac
\def\R{\mathbb R}
\def\P{\mathbb P}
\def\E{\mathbb E}

\def\eqd{\buildrel\hbox{\footnotesize d}\over =}
\def\tod{\buildrel\hbox{\footnotesize d}\over \longrightarrow}

\newcommand{\g}{\bar{g}}
\usepackage{color}
\def\1{\mathbf{1}}
\begin{document}
	\title[Principal minors ]{Principal minors of Gaussian orthogonal ensemble}
	
	\author{Renjie Feng}
	\author{Gang Tian}
	\author{Dongyi Wei}
	\author{Dong Yao}
	
	\address{Sydney Mathematical Research Institute, The University of Sydney, Sydney, Australia, 2006. }
	\email{renjie.feng@sydney.edu.au}
	\address{Beijing International Center for Mathematical Research and School of Mathematical Sciences, Peking University, Beijing, China, 100871.}
	\email{gtian@math.pku.edu.cn}
	\address{School of Mathematical Sciences, Peking University, Beijing, China, 100871.}
	\email{jnwdyi@pku.edu.cn}
	\address{Research Institute of Mathematical Science, Jiangsu Normal University, Xuzhou, China, 221116. }
	\email{wonderspiritfall@gmail.com}

	\maketitle
	\begin{abstract} In this paper,  we study the extremal process of the maxima of all the largest eigenvalues of principal
		minors  of the classical Gaussian orthogonal ensemble (GOE). We prove that the fluctuation of the maxima  is given by the Gumbel distribution in the limit. We also derive the limiting joint distribution of the maxima and the corresponding  eigenvector, which implies that these two random variables are asymptotically independent. 
		
	\end{abstract}
	\section{Introduction}\label{intro}

	Motivated by   high-dimensional statistics and signal processing, the authors in \cite{CJL}  derived
	the growth order of the extremal process of the maxima of  all the largest eigenvalues of the principal minors of the classical random matrices of GOE and Wishart matrices, where the results have   applications for the construction of the compressed  sensing matrices as in \cite{TC2}. In this paper, we further study the fluctuation of such maxima  for the GOE case. Our main result is that the fluctuation is given by the Gumbel distribution with some  Poisson structure involved in the limit, and we also derive the limiting joint distribution of such maxima and its corresponding eigenvector which indicates that these two random variables are asymptotically independent.  
	
	The Gaussian orthogonal ensemble (GOE) is the Gaussian measure defined on the space of real symmetric matrices, i.e., $G=({g}_{ij})_{1\leq i,j\leq n} $ is a
	symmetric matrix whose upper triangular entries are independent real Gaussian variables with
	the following distribution\begin{align*}& {g}_{ij}\eqd\left\{\begin{array}{l}N_{\mathbb R}(0,2)\ \text{if}\ i=j;\\ N_{\mathbb R}(0,1)\ \text{if}\ i<j.\end{array}\right.
	\end{align*} 
	Let $\lambda_1(G)>\lambda_2(G)>\cdots >\lambda_n(G)$ be eigenvalues of GOE, then the distribution of these eigenvalues is invariant under the orthogonal group action and the joint density  is 
	\begin{equation}\label{jointgoe} \frac{1}{Z_n} \prod_{k=1}^n e^{-\frac{ 1}{4}\lambda_k^2}\prod_{i<j}\big(\lambda_i-\lambda_j\big),  \end{equation}
	where \begin{equation}\label{normgoe}Z_n=2^{n(n+1)/4}(2\pi)^{n/2}(n!)^{-1}\prod_{j=1}^{n}\frac{\Gamma(1+j/2)}{\Gamma(3/2)}\end{equation} is the partition function. And the limit of the empirical measure of these eigenvalues is given by the classical semicircle law \cite{AGZ}. 
	
	Let's first introduce some notations in order to present our main results. Given symmetric matrices $G=({g_{ij}})_{1\leq i\leq j \leq n}$ sampled from GOE, for $ \alpha\subset\{1,\cdots,n\}$ with cardinality $| \alpha|=m\in \mathbb Z_+$, we denote ${G}_ \alpha=({g}_{ij})_{i,j\in \alpha}$ as the principal minors  of $G$ of size $m\times m$, then $G_ \alpha$ is also  symmetric.  Let $\lambda_1(G_ \alpha)>\lambda_2(G_ \alpha)>\cdots>\lambda_m(G_ \alpha)$ be  the ordered eigenvalues of $G_\alpha$. Now we define the extremal process of  the maxima of all the largest eigenvalues of the principal minors as \begin{align*}&  {T}_{m,n}=\max_{ \alpha\subset\{1,\cdots,n\},| \alpha|=m}\lambda_1( {G}_ \alpha). 
	\end{align*} 
	In \cite{CJL},  the authors studied the asymptotic properties of ${T}_{ m,n}$ and proved that under the assumption that $m$ fixed or $m\to+\infty$ with $m=o\Big(\frac{(\ln n)^{1/3}}{\ln\ln n}\Big)$, it holds
	$$T_{m, n}-2\sqrt{m\ln n} \to 0$$
	in probability as $n\to + \infty$. 
	


	
	In this paper, we further derive the fluctuation of ${T}_{m,n}$ when $m$ is fixed as $n$ tends to infinity. Our first result is the following 
	\begin{thm}\label{main} For \text{GOE},  we have the following convergence in distribution \begin{equation}\label{tmnlim} {T}_{m,n}^2-4m\ln n-2(m-2)\ln\ln n\tod Y
		\end{equation}  as $n\to+\infty$ for fixed $m$, where the random
		variable $Y$ has the Gumbel distribution function\begin{align*}& F_{Y}(y)=\exp(-c_{m}e^{-y/4}),\quad y\in\R.
		\end{align*} Here, the constant $$c_m=\frac{(2m)^{(m-2)/2}K_{m}}{(m-1)!2^{3/2}\Gamma(1+m/2)},$$
		where $K_m=\mu(\mathcal S_m)$ is the probability of the event \begin{equation}\label{hm}\mathcal S_m:=\Big\{x\in S^{m-1}: \sum_{j\in\beta}x_j^2\leq\sqrt{k/m},\,\,\forall\, \beta \subset \{1,\cdots, m\}\,\,\text{with }\,\, \forall \,1\leq |\beta|=k<m   \Big\}\end{equation} under the uniform distribution  $\mu$ on the unit sphere $S^{m-1}$. 
		In particular, $$c_1=\frac 1{2\sqrt \pi},\quad  c_2=-\frac1{2\sqrt 2}+\frac {\sqrt 2}{\pi} \arcsin  \left(\big(\frac 12\big)^{1/4}\right).$$

	\end{thm}

	The following classical result in {the theory of extremal processes} is regarded as a special case of our main result.  Let $a_1,..., a_n$ be i.i.d. Gaussian random variables $N_{\mathbb R}(0,2)$, for the extremal process 
	$$M_n:=\max\{a_1,...,a_n\},$$ let $$a_{n}=2\sqrt{ \ln n} $$
	and $$b_{n} = 2\sqrt{ \ln n}-\frac{\ln \ln n+\ln (4 \pi)}{2 \sqrt{ \ln n}}.$$
	Then for any $y \in \mathbb{R}$,   the following classical result holds (Theorem 1.5.3 in \cite{LLr})
	$$
	\lim _{n \to +\infty} \mathbb{P}\left[a_{n}\left(M_{n}-b_{n}\right) \leq y\right]=e^{-e^{-y/2}}.
	$$
	One can check that Theorem \ref{main} for $T^2_{1, n}$ when $m=1$ is   equivalent to this classical result for $M_n$.  In this sense, our result is fundamental which  is  a natural generalization of such  classical result for the extremal process of the scalar-valued random variables to the matrix-valued random variables with correlations. 
	
	The constant $K_m$ appearing in the extremal index $c_m$ represents the strength of the correlation between different principal minors. For a comprehensive overview of the extremal theory of stochastic processes, we refer to \cite{LLr, LR}.	
	
	
	Our proof for Theorem \ref{main} also implies the joint distribution of the maxima of the largest  eigenvalues of principal minors and its corresponding eigenvector. To be more precise, 
	given $n$,  let $v^*\in S^{m-1}$ be the unit eigenvector corresponding to the largest eigenvalue of the principal minor that attains the maxima  $T_{m,n}$. By symmetry, $-v^*$ is also the corresponding eigenvector. Now we have the  following limit for the joint distribution of $(T_{m,n},  v^*)$.   
	\begin{thm}\label{main2}  Given any $y\in \R$ and symmetric Borel set $Q\subset S^{m-1}$ such that $-Q=Q$, let $y_m^2=4m\ln n+2(m-2)\ln\ln n+y$, 
		then the joint distribution satisfies 
		\begin{equation}\label{value+vec}
			\P(T_{m,n}>y_m,  v^* \in Q ) \to  (1-F_Y(y))\nu (Q)
		\end{equation}
		as $n \to+\infty$, which implies that $T_{m,n}$  and $v^*$ are asymptotically independent. 
		Here, $\nu$ is the uniform distribution on    $\mathcal S_m$, i.e., $$\nu(Q)=\mu(Q\cap \mathcal S_m)/\mu(\mathcal S_m). $$
		
	\end{thm}
	
	The study of principal minors of random matrices holds significant importance and finds applications in various fields such as statistics, machine learning, deep learning, information theory, and more.  In the following, we only list two of them. 
	
	In {compressed sensing}, one has to recover an input vector $f$ from the corrupted measurements $y=A f+e$, where $A$ is a coding matrix and $e$ is an arbitrary and unknown vector of errors. The  famous result by Cand\`{e}s-Tao \cite{TC2} is that if  the coding matrix $A$ satisfies 
	the restricted isometry property  (Definition 1.1 in \cite{TC2}),  then the input $f$ is the unique solution to some $\ell_{1}$-minimization problem 
	provided that the support $S$ (the number of nonzero entries) of errors $e$ is not too large. 
	Therefore, one of the major goals in compressed sensing is to construct the coding
	matrix $A$ that satisfies the restricted isometry property.  In Section 3 of \cite{TC2},  Cand\`{e}s-Tao proved that the Gaussian random matrices $A$ can satisfy such property with overwhelming probability,  and they can derive the estimate about the support $S$ 
	via the probabilistic estimate on the maxima of the largest eigenvalues of principal minors of  Wishart matrices $A^TA$. 
	A simple proof based on the concentration measure theory is present in \cite{B}.

	In statistics, the problem of finding the largest eigenvalue over  principal minors  (of some fixed order)  of the sample covariance matrices is also called \emph{sparse principal component analysis}.   
	Gamarnik and Li \cite[Section 6]{GL} asked about the distribution of the largest eigenvalue as well as the  algorithmic complexity of this problem, i.e., whether there exists efficient algorithms to find the principal minor which attains the (near-)optimal value for the largest eigenvalue. Consider   the Wishart version of this problem, i.e.,
	let $X$ be an $n\times p$ matrices with independent $N_{\mathbb R}(0,1)$ distributed entries and take $R_{m,n,p}$ to be the maximal eigenvalue of all $m$-th order principal minors of $X^TX$. It is proved in \cite{CJL}	that 
	$$
	\frac{R_{m,n,p}-n}{\sqrt{n}}-\sqrt{2m\log p} \to 0
	$$
	in probability as $n\to\infty$, under  conditions requiring $\log p=o(n^{1/2})$. 
	It seems that our method in the current paper can also be applied to study the Wishart case and it's expected that some Gumbel  fluctuation  will be observed as well. We will postpone this to  further investigations.

	\bigskip

	\emph{Notation.} In this paper, $c$, $C$ and $C'$ stand for  positive constants, but their values may change from line to line. For simplicity, the notation $a_n\sim b_n$ means $\lim_{n\to+\infty} a_n/b_n=1$.  \section{Proof of Theorem \ref{main}} \label{main1}
	In this section,   we will prove   Theorem \ref{main}   by assuming some technical lemmas where the proofs of these lemmas are postponed to \S \ref{3}.

	The proof of Theorem \ref{main} is based on  Lemma \ref{lem:Poisson} with the proof given in  \cite{AGG} by the Stein-Chen method. It provides a   criteria to prove the convergence of the total number of   occurrences of the point process to the Poisson distribution,  and thus it provides a method to derive the distribution for some extremal processes. \begin{lem}\label{lem:Poisson}
		Let $I$ be an index set, and for $\alpha\in I $, let $ X_{\alpha}$ be a Bernoulli random variable with $p_{\alpha}=\mathbb{P}(X_{\alpha}=1)=1-\mathbb{P}(X_{\alpha}=0)$. For each $ \alpha\in I$, let $N_{\alpha}$ be a subset of $I$ with $ \alpha\in N_{\alpha}$, that is, $\alpha\in N_{\alpha}\subset I$. Let $$S=\sum_{\alpha\in I}X_{\alpha},\,\, \lambda=\mathbb E S=\sum_{\alpha\in I}p_{\alpha}\in (0, +\infty),$$ 
		let $Z$ be the Poisson random variable with intensity $\mathbb E Z=\lambda$, 
		then it holds that 
		\begin{align*}\|\mathcal L(S)-\mathcal L (Z)\|\leq 2(b_1+b_2+b_3),
		\end{align*} and the probability of no occurrence has the estimate
		\begin{align*}&|\mathbb{P}(S=0)-e^{-\lambda}|\\=&|\mathbb{P}(X_{\alpha}=0,\ \forall\ \alpha\in I)-e^{-\lambda}|\\ \leq &  \min(1,\lambda^{-1})(b_1+b_2+b_3).
		\end{align*} Here, $\|\mathcal L(S)-\mathcal L (Z)\|$ is the total variation distance between the distributions $S$ and $Z$, 
		and \begin{align*}&b_1=\sum_{\alpha\in I}\sum_{\beta\in N_{\alpha}}p_{\alpha}p_{\beta},\\& b_2=\sum_{\alpha\in I}\sum_{\alpha\neq\beta\in N_{\alpha}}\mathbb{E}[X_{\alpha}X_{\beta}],\\&b_3=\sum_{\alpha\in I}\mathbb{E}|\mathbb{E}[X_{\alpha}|\sigma(X_{\beta},\beta\not\in N_{\alpha})]-p_{\alpha}|,
		\end{align*}and $\sigma(X_{\beta},\beta\not\in N_{\alpha})$ is the $ \sigma$-algebra generated by $\{X_{\beta},\beta\not\in N_{\alpha}\}.$ In particular, if $X_{\alpha} $ is
		independent of $\{X_{\beta},\beta\not\in N_{\alpha}\}$ for each $\alpha$, then $b_3=0$.
	\end{lem}
	
	One may think of $N_\alpha$ as  a `neighborhood of dependence' for $\alpha$ such that $X_\alpha$
	is independent or nearly independent of all of $X_\beta$ for $\beta\notin N_\alpha$. And Lemma \ref{lem:Poisson} indicates that when $b_1$, $b_2$ and $b_3$ are all small enough, then $S$ which is the total number of occurrences    tends to a Poisson distribution. 

	\subsection{Proof of Theorem \ref{main}}
	
	
	Given any $\ell\times \ell$ symmetric matrix $S$, we rearrange  the eigenvalues of $S$ in descending order $ \lambda_1(S)\geq\cdots\geq\lambda_\ell(S)$, and we denote \begin{equation}\label{trace2}|S|^2:=\Tr S^2\end{equation} and \begin{equation}\label{trace12}\lambda_1^*(S):=[(\lambda_1(S)^2+|S|^2)/2]^{\f12}.\end{equation}
	For fixed $m\geq 2$, we define  the index set $$I_m=\{\alpha \subset\{1,\cdots,n\},|\alpha|=m\}$$ and   the neighborhood set $$N_{\alpha}=\{\beta\in I_m:\alpha\cap\beta\neq\emptyset\}\,\,\,\mbox{for}\,\,\, \alpha\in I_m.$$ 
	Throughout the article, for any fixed real number $y$, 
	we define $$y_k^2:=4k\ln n+2k\ln\ln n,\ y_k>0,\ 1\leq k<m$$ and 
	$$y_m^2:=4m\ln n+2(m-2)\ln\ln n+y,\ y_m>0.$$
	For  symmetric matrices $G=({g_{ij}})_{1\leq i\leq j \leq n}$ sampled from GOE, for $\alpha\in I_m$ where $|\alpha|=m$, we denote ${G}_\alpha=({g}_{ij})_{i,j\in \alpha}$ as the principal minor  of size $m\times m$ and we define the event \begin{align}\label{defAS}
		A_{\alpha}=\Big\{\lambda_1({G}_\alpha)>y_{m};\ \ \lambda_1^*({G}_{\beta})\leq y_{k},\ \forall\ 1\leq k<m,\ \beta\subset \alpha,\ |\beta|=k\Big\}.
	\end{align}Recall the definition of ${T}_{m,n}$ in \S \ref{intro},   we first have 
	\begin{equation}\label{step1}
		\begin{split}0&\leq \mathbb{P}(\cap_{\alpha\in I_m}A_{\alpha}^c)-\mathbb{P}({T}_{m,n}\leq y_{m})\\&  \leq \mathbb{P}(\lambda_1^*({G}_{\beta})> y_{k}\ \text{for some}\ \beta\in I_k,\ 1\leq k<m)\\&\leq \sum_{k=1}^{m-1}\sum_{\beta\in I_k}\mathbb{P}(\lambda_1^*({G}_{\beta})> y_{k})\\& =\sum_{k=1}^{m-1}{n\choose k}\mathbb{P}(\lambda_1^*({G}_{\{1,\cdots,k\}})> y_{k}).
		\end{split}
	\end{equation}
	Now we need the following lemma and we postpone its proof to the next section. 
	\begin{lem}\label{lem:WS}
		For fixed $k\geq $1, there exists a constant $C>0$ (depending on $k$) so that for all  $x>1$, \begin{align}&\mathbb{P}(|{G}_{\{1,\cdots,k\}}|^2> x^2)\leq Cx^{k(k+1)/2-2}e^{-x^2/4}, \label{w1}\\&\mathbb{P}(\lambda_1({G}_{\{1,\cdots,k\}})> x)\leq Cx^{k-2}e^{-x^2/4},\label{w2}\\&\mathbb{P}(\lambda_1^*({G}_{\{1,\cdots,k\}})> x)\leq Cx^{k-2}e^{-x^2/4} \label{w3}. 
		\end{align}
	\end{lem}
	Using \eqref{w3} we have  \begin{equation}\label{limit0}\sum_{k=1}^{m-1}{n\choose k}\mathbb{P}(\lambda_1^*({G}_{\{1,\cdots,k\}})> y_{k})\leq C\sum_{k=1}^{m-1}n^ky_k^{k-2}e^{-y_k^2/4}\leq C/\ln n.
	\end{equation}
	Combining this with \eqref{step1} we get
	\begin{equation}\label{ooo}0\leq \mathbb{P}(\cap_{\alpha\in I_m}A_{\alpha}^c)-\mathbb{P}({T}_{m,n}\leq y_{m})\leq C/\ln n, 
	\end{equation}
	and thus we have 
	\begin{equation}\label{lim2}\lim_{n\to+\infty} \mathbb{P}({T}_{m,n}\leq y_m)=\lim_{n\to+\infty} \mathbb{P}(\cap_{\alpha\in I_m}A_{\alpha}^c).
	\end{equation}
	Therefore, it's enough to derive the limit of $\mathbb{P}(\cap_{\alpha\in I_m}A_{\alpha}^c)$ to prove Theorem \ref{main}. 
	By Lemma \ref{lem:Poisson}, we have \begin{equation}\label{dd} |\mathbb{P}(\cap_{\alpha\in I_m}A_{\alpha}^c)-e^{-t_n}|\leq b_{n,1}+b_{n,2},
	\end{equation}
	where
	\begin{align*}& t_n={n\choose m}\mathbb{P}(A_{\{1,\cdots,m\}}),\\&  b_{n,1}\leq {n\choose m}\left({n\choose m}-{n-m\choose m}\right)\mathbb{P}(A_{\{1,\cdots,m\}})^2,\\& b_{n,2}\leq \sum_{k=1}^{m-1}{n\choose k}{n-k\choose m-k}{n-m\choose m-k}\mathbb{P}(A_{\{1,\cdots,k, \cdots, m\}}\cap A_{\{1,\cdots,k,m+1,\cdots,2m-k\}}).
	\end{align*}
	Using \eqref{w2} we have
	\begin{align*}&\mathbb{P}(A_{\{1,\cdots,m\}})\leq \mathbb{P}(\lambda_1({G}_{\{1,\cdots,m\}})> y_{m})\leq Cy_m^{m-2}e^{-y_m^2/4}\leq Cn^{-m}.
	\end{align*} And thus we have 
	\begin{align*}&b_{n,1}\leq Cn^{2m-1}\mathbb{P}(A_{\{1,\cdots,m\}})^2\leq C'n^{-1}, 
	\end{align*} which tends to $0$ in the limit.

	It remains to find the limit of $t_n$ and show that $b_{n,2}$ tends to $0$ in order to complete the proof of Theorem \ref{main}. 
	
	Let $ \alpha=\{1,\cdots,m\}, \gamma=\{m-k+1,\cdots,m\}, \zeta=\{m-k+1,\cdots,2m-k\}$, then $\alpha\cap\zeta=\gamma$, $|\alpha|=|\zeta|=m$ and $|\gamma|=k$. By rearranging the indices,  we have 
	\begin{align*}&\mathbb{P}(A_{\{1,\cdots,k, \cdots, m\}}\cap A_{\{1,\cdots,k,m+1,\cdots,2m-k\}})\\=& \mathbb{P}(A_{\{1,\cdots,m\}}\cap A_{\{m-k+1,\cdots,2m-k\}})\\ \leq& \mathbb{P}(\lambda_1({G}_{\alpha})>y_{m},\ \lambda_1( {G}_{\zeta})>y_{m},\\& \lambda_1^*( {G}_{\alpha\setminus\gamma})\leq y_{m-k},\ \lambda_1^*( {G}_{\zeta\setminus\gamma})\leq y_{m-k},\ \lambda_1^*( {G}_{\gamma})\leq y_{k})\\ =& \mathbb{E}[\mathbb{P}(\lambda_1( {G}_{\alpha})>y_{m},\ \lambda_1( {G}_{\zeta})>y_{m},\ \lambda_1^*( {G}_{\alpha\setminus\gamma})\leq y_{m-k},\ \lambda_1^*( {G}_{\zeta\setminus\gamma})\leq y_{m-k}| {G}_{\gamma} )\\&\times\mathbf{1}_{\{\lambda_1^*( {G}_{\gamma})\leq y_{k}\}}]\\ =& \mathbb{E}[\mathbb{P}(\lambda_1( {G}_{\alpha})>y_{m},\ \lambda_1^*( {G}_{\alpha\setminus\gamma})\leq y_{m-k}| {G}_{\gamma} )^2\mathbf{1}_{\{\lambda_1^*( {G}_{\gamma})\leq y_{k}\}}].
	\end{align*} 
	The following lemma will imply that $b_{n,2}$ tends to $0$ as $n \to+\infty$. 
	\begin{lem}\label{lem:WS1}
		For $\alpha\in I_m,\ \gamma\subset\alpha,\ |\gamma|=k,\ 1\leq k<m,\ \beta=\alpha\setminus\gamma,\ x>1,\ \delta,\delta' \in(0,1)$, 
		then there are some constants $C$ and $C'$ (depending on $m$, $\delta$ and $\delta'$) such that
		\begin{align}
			\label{alpha} &\mathbb{P}(\lambda_1({G}_{\alpha})> x|{G}_{\beta}, {G}_{\gamma})\mathbf{1}_{\{\lambda_1^*(G_{\beta})\leq (1-\delta)x,\lambda_1^*(G_{\gamma})\leq (1-\delta')x\}} \\
			\leq& Cx^{-1}(x/(\lambda_1^*(G_{\beta})+1)+x/(\lambda_1^*(G_{\gamma})+1))^{k(m-k)-1}e^{-(x-\lambda_1^*(G_{\beta}))(x-\lambda_1^*(G_{\gamma}))/2}\nonumber \end{align} 
		and 
		\begin{align}
			&\mathbb{P}(\lambda_1({G}_{\alpha})> x,\ \lambda_1^*({G}_{\beta})\leq (1-\delta)x| {G}_{\gamma})\mathbf{1}_{\{\lambda_1^*(G_{\gamma})\leq (1-\delta')x\}}\label{alphabeta}\\ \leq& C'x^{m-k-2}(x/(\lambda_1^*(G_{\gamma})+1))^{k(m-k)-1}e^{((\lambda_1^*(G_{\gamma}))^2-x^2)/4} \nonumber .
		\end{align}
	\end{lem} 
	By assuming Lemma \ref{lem:WS1}, we have\begin{align}& \mathbb{E}[\mathbb{P}(\lambda_1( {G}_{\alpha})>y_{m},\ \lambda_1^*( {G}_{\alpha\setminus\gamma})\leq y_{m-k}| {G}_{\gamma} )^2\mathbf{1}_{\{\lambda_1^*( {G}_{\gamma})\leq y_{k}\}}] \nonumber\\ \leq &C\mathbb{E}[y_m^{2m-2k-4}(y_m/(\lambda_1^*(G_{\gamma})+1))^{2k(m-k)-2}e^{((\lambda_1^*(G_{\gamma}))^2-y_m^2)/2}\mathbf{1}_{\{\lambda_1^*( {G}_{\gamma})\leq y_{k}\}}]\nonumber\\ =&:C \E[f(\lambda_1^*(G_{\gamma})))\mathbf{1}_{\{\lambda_1^*( {G}_{\gamma})\leq y_{k}\}}] \label{square1}
	\end{align}
	where we define 
	\begin{equation}\label{flem3}
		f(t)=y_m^{2m-2k-4}(y_m/(t+1))^{2k(m-k)-2}e^{(t^2-y_m^2)/2},\,\,\, t\geq 0.
	\end{equation}
	Integration by parts, we have
	\begin{equation}\label{ef2}
		\begin{split}
			\E(f(\lambda_1^*(G_{\gamma}))\1_{\lambda_1^*(G_{\gamma})\leq y_k})=&\int_0^{y_k} f'(t)\P(\lambda_1^*(G_{\gamma})>t)dt\\
			&-f(y_k)\P(\lambda_1^*(G_{\gamma})>y_k) +f(0).
		\end{split}
	\end{equation}
	Note that 
	\begin{equation}\label{f'tbd}
		\begin{split}
			f'(t)&=-(2k(m-k)-2)(t+1)^{-1}f(t)+tf(t)\leq tf(t).\\
		\end{split}
	\end{equation}
	Hence, using \eqref{w3}  for $t>1$ and \eqref{f'tbd},   we  have 
	\begin{equation}\label{intf'bd}
		\begin{split}
			&\int_0^{y_k} f'(t)\P(\lambda_1^*(G_{\gamma})>t)dt\\
			\leq &\int_0^1 tf(t) dt +\int_1^{y_k}f'(t)\P(\lambda_1^*(G_{\gamma})>t)dt\\
			\leq&\max_{0\leq t\leq 1}f(t)+C\int_1^{y_k} t y_m^{2m-2k-4}(y_m/(t+1))^{2k(m-k)-2}e^{(t^2-y_m^2)/2}t^{k-2} e^{-t^2/4}dt.
		\end{split}
	\end{equation}
	Therefore, by \eqref{ef2} we further have \begin{equation}\label{ef}
		\begin{split}
			&\E(f(\lambda_1^*(G_{\gamma}))\1_{\lambda_1^*(G_{\gamma})\leq y_k}) \leq 2\max_{0\leq t\leq 1}f(t) \\
			&+C\int_1^{y_k} t y_m^{2m-2k-4}(y_m/(t+1))^{2k(m-k)-2}e^{(t^2-y_m^2)/2}t^{k-2} e^{-t^2/4}dt.
		\end{split}
	\end{equation}
	We now separate this integration into $1<t<y_k/2$ and $y_k/2<t<y_k$. For $1<t<y_k/2$ the integrand is bounded by
	$$
	C y_m^{2m-k-6} y_m^{2k(m-k)-1} e^{y_k^2/16-y_m^2/2},
	$$
	where we used the fact that $y_k\leq y_m$ for $n$ large enough.
	For $y_k/2<t<y_k$, we can bound the integrand by 
	$$
	C y_m^{2m-k-6}  te^{t^2/4-y_m^2/2},
	$$
	where we used the fact that $y_m/y_k \leq 2 \sqrt{m/k}$ as $n$ large enough. 
	Therefore,  as $n$ large enough, we have 
	\begin{equation}\label{0toy_k}
		\begin{split}
			&\int_1^{y_k} f'(t)\P(\lambda_1^*(G_{\beta})>t)dt\\
			\leq & Cy_m^{2m-k-6}e^{-y_m^2/2}\left( \int_1^{y_k/2}  y_m^{2k(m-k)-1} 
			e^{y_k^2/16} dt + \int_{y_k/2}^{y_k} t e^{t^2/4}dt  \right)\\
			\leq & C y_m^{2m-k-6}e^{-y_m^2/2}
			\left(  y_m^{2k(m-k)}  e^{y_k^2/16}+e^{y_k^2/4}
			\right)
			\\
			\leq &Cy_m^{2m-k-6}e^{y_k^2/4-y_m^2/2},
		\end{split}
	\end{equation}  where in the last inequality we used the fact that $ y_m^{2k(m-k)}  e^{-3y_k^2/16}$ can be  bounded from above uniformly for all $n$. 
	
	The definition of $f(t)$ in \eqref{flem3} and the fact that $y_m/y_k \sim \sqrt{m/k}$ imply
	\begin{equation}\label{f0}
		\max_{0\leq t\leq 1} f(t)\leq y_m^{2m-2k-4} y_m^{2k(m-k)-2}e^{(1-y_m^2)/2}\leq 
		Cy_m^{2m-k-6}e^{y_k^2/4-y_m^2/2}
	\end{equation} as $n$ large enough. 
	It follows from \eqref{square1}, \eqref{ef}, \eqref{0toy_k} and \eqref{f0} that
	\begin{equation}\label{final1}
		\begin{split}
			&\mathbb{E}[\mathbb{P}(\lambda_1( {G}_{\alpha})>y_{m},\ \lambda_1^*( {G}_{\alpha\setminus\gamma})\leq y_{m-k}| {G}_{\gamma} )^2\mathbf{1}_{\{\lambda_1^*( {G}_{\gamma})\leq y_{k}\}}] \\ 
			\leq &  Cy_m^{2m-k-6}e^{y_k^2/4-y_m^2/2}
		\end{split}
	\end{equation}as $n$ large enough. 
	Therefore, we have 
	\begin{equation}\label{bn2ctl}
		\begin{split}
			b_{n,2}& \leq  \sum_{k=1}^{m-1}{n\choose k}{n-k\choose m-k}{n-m\choose m-k}\mathbb{P}(A_{\{1,\cdots,k, \cdots, m\}}\cap A_{\{1,\cdots,k,m+1,\cdots,2m-k\}})\\&\leq \sum_{k=1}^{m-1}n^{2m-k}\mathbb{P}(A_{\{1,\cdots,m\}}\cap A_{\{1,\cdots,k,m+1,\cdots,2m-k\}})\\&\leq C\sum_{k=1}^{m-1}n^{2m-k}y_m^{2m-k-6}e^{y_k^2/4-y_m^2/2}\\&\leq C\sum_{k=1}^{m-1}y_m^{2m-k-6}(\ln n)^{k/2-m+2}\leq C(\ln n)^{-1},
		\end{split}
	\end{equation} which tends to $0$ as $n\to+\infty$.
	
	The following lemma gives the limit of $ t_n$. 
	\begin{lem}\label{t_nlimit}
		We have the limit 
		\begin{equation}
			\lim_{n\to +\infty}t_n = c_me^{-y/4},
		\end{equation}
		where $$c_m=\frac{(2m)^{(m-2)/2}K_{m}}{(m-1)!2^{3/2}\Gamma(1+m/2)}, $$
		where $K_m$ is  the constant defined in Theorem \ref{main}. 
	\end{lem}
	By assuming Lemma \ref{t_nlimit},  by \eqref{dd} together with the facts that $b_{n,1}\to0$ and $b_{n,2}\to 0$, we can conclude that $$\lim_{n\to+\infty} \mathbb{P}(\cap_{\alpha\in I}A_{\alpha}^c)= \exp(-c_me^{-y/4}).$$
	By \eqref{lim2}, this further implies $$\lim_{n\to+\infty} \mathbb P(T_{m,n }\leq y_m)=\exp(-c_me^{-y/4}),$$ 
	which proves Theorem \ref{main}.


	\section{Proofs of lemmas} \label{3}
	In this section, we will prove Lemma \ref{lem:WS}, Lemma \ref{lem:WS1} and Lemma \ref{t_nlimit}, and thus we complete the proof of Theorem \ref{main}.

	\subsection{Proof of Lemma \ref{lem:WS}}
	\begin{proof}
		We simply have the following estimates.  For $s\in\mathbb R$,  there are some constants $C, C'>0$ depending on $s$ such that for all $x>1$ we have
		\begin{equation}
			\int_x^{+\infty} r^s\exp(-r)dr\leq C x^s\exp(-x) \label{r} \end{equation} and 
		\begin{equation} \int_x^{+\infty} r^s \exp(-r^2/2)dr\leq C' x^{s-1}\exp(-x^2/2). \label{r^2}
		\end{equation}
		Now let $g_1,\ldots, g_\ell $ be $\ell $ independent $N_{\mathbb{R}}(0,1)$ random variables, then for all $t\geq 1$, we have 
		\begin{equation} \label{xi^2}
			\P\Big(\sum_{i=1}^\ell g_i^2 \geq t \Big) \leq Ct^{\ell /2-1}\exp(-t/2),
		\end{equation} where $C>0$ only depends on $\ell$. 
		The proof of \eqref{xi^2} follows if we combine the estimate \eqref{r} and the fact that the probability density of the chi-squared distribution $\chi^2(\ell):=\sum_{i=1}^\ell  g_i^2$ with $\ell$ degrees of freedom is given by $$\frac 1{2^{\ell /2}\Gamma(\ell /2)}x^{\ell /2-1}e^{-x/2}.$$
		Lemma \ref{lem:WS} holds obviously when $k=1$ (note that $g_{ii}\eqd N_{\mathbb R}(0,2)$), and thus in the followings we consider the case when $k\geq 2$. By the definition of the principal minor, $G_{1,\ldots, k}=(g_{ij})_{1\leq i,j\leq k}$ is also sampled from GOE. Now let $\g_{ij}=g_{ij}$ if $j\neq i$ and $g_{ij}/\sqrt 2$ otherwise, then $\g_{ij},1\leq i\leq j\leq k$ are i.i.d. $N_\mathbb{R}(0,1)$ random variables. 
		To prove \eqref{w1}, we note that by definition,
		\begin{align*}
			|G_{1,\ldots, k}|^2&=\Tr(G_{1,\ldots, k}^2)=\sum_{i,j=1}^k g_{ij}^2=2\sum_{1\leq i\leq j\leq k} \g_{ij}^2,
		\end{align*}
		therefore,  by \eqref{r} and \eqref{xi^2} we have 
		\begin{align*}
			\P(|G_{1,\ldots, k}|^2>x^2)&=\P\Big(\sum_{1\leq i\leq j\leq k} \g_{ij}^2> x^2/2\Big) 
			\leq Cx^{k(k+1)/2-2}e^{-x^2/4}. 
		\end{align*}
		This proves \eqref{w1}. 
		To prove \eqref{w2},  by formula \eqref{jointgoe}, the joint density
		of eigenvalues $\lambda_1\geq  \cdots \geq  \lambda_k$ of ${G}_{\{1,\cdots,k\}}$ is $$
		J_k(\lambda_1, \ldots,  \lambda_k)=\frac{1}{Z_k} \prod_{1\leq i<j\leq k}(\lambda_i-\lambda_j)\exp\left(-
		\frac{\sum_{i=1}^k\lambda_i^2}{4}\right).
		$$
		Note that if $\lambda_1>x>1$, we have
		$
		\lambda_1-\lambda_j \leq (\lambda_1+1)(\abs{\lambda_j}+1)\leq 
		2\lambda_1 (\abs{\lambda_j}+1), 
		$
		and thus we have 
		\begin{equation}\label{updensity}
			\prod_{1\leq i<j\leq k}(\lambda_i-\lambda_j)\leq C \lambda_1^{k-1} 
			\prod_{2\leq i\leq k} (\abs{\lambda_i}+1) \prod_{2\leq i<j\leq k}(\lambda_i-\lambda_j), 
		\end{equation}
		which further implies 
		\begin{align*} & \int_{\lambda_1>x}J_k(\lambda_1, \ldots, \lambda_k)d\lambda_1\cdots d\lambda_k \leq  C\int_{x}^{\infty}\lambda_1^{k-1} \exp(-\lambda_1^2/4)d\lambda_1 \\&   \times \int_{+\infty>\lambda_2>\cdots>\lambda_k>-\infty} \prod_{2\leq i\leq k} (\abs{\lambda_i}+1) \prod_{2\leq i<j\leq k}(\lambda_i-\lambda_j)\exp\left(-
			\frac{\sum_{i=2}^k\lambda_i^2}{4}\right)d\lambda_2\cdots d\lambda_k, 
		\end{align*}
		which is bounded from above by $C'x^{k-2}\exp(-x^2/4)$ by \eqref{r^2}, thereby proving \eqref{w2}. 
		
		Now we prove \eqref{w3}. 
		By definition of $\lambda_1^*(G_{1,\ldots, k} )$, we have 
		\begin{align*}
			&\P(\lambda_1^*(G_{1,\ldots, k} )>x) \\=&\P\left(\sum_{i=2}^k \lambda_i^2+2\lambda_1^2>2x^2\right)\\
			\leq &\sum_{y=0}^{\lfloor \sqrt{2}x\rfloor+1} \P\left(2\lambda_1^2>(2x^2-(y+1)^2)_+, \,\,y^2\leq \sum_{i=2}^k \lambda_i^2\leq (y+1)^2\right)\\
			&+\P\left(\sum_{i=2}^k \lambda_i^2>2x^2\right)\\
			:=&I_1+I_2.
		\end{align*}
		By  the fact $
		\lambda_1-\lambda_j \leq (|\lambda_1|+1)(\abs{\lambda_j}+1)$ for all $2\leq j\leq k$,  we first have 
		\begin{equation}\label{stis} \prod_{1\leq i<j\leq k}(\lambda_i-\lambda_j)\leq C (|\lambda_1|+1)^{k-1} 
			\prod_{2\leq i\leq k} (\abs{\lambda_i}+1)  \prod_{2\leq i<j\leq k}(\lambda_i-\lambda_j). \end{equation}
		Then by the inequality of arithmetic and geometric means, for $\sum_{i=2}^k \lambda_k^2\geq 1$,  
		we have 
		\begin{equation}\label{sdsds}\begin{split}
				&\prod_{2\leq i\leq k} (\abs{\lambda_i}+1) \prod_{2\leq i<j\leq k}(\lambda_i-\lambda_j) \leq  \prod_{2\leq i\leq k} (\abs{\lambda_i}+1)\prod_{2\leq i<j\leq k}(|\lambda_i|+|\lambda_j|)\\ \leq & \left(\frac{ 1+\sum_{i=2}^k |\lambda_i|}{k/2} \right)^{k(k-1)/2}\leq   \left(\frac{ 1+\sqrt{k-1}\sqrt{\sum_{i=2}^k \lambda_i^2}}{k/2} \right)^{k(k-1)/2}\leq   C  \left(\sum_{i=2}^k \lambda_i^2 \right)^{k(k-1)/4}. \end{split}
		\end{equation}
		Therefore, for $x>1$,  combining \eqref{stis} and \eqref{sdsds},  we can bound $I_2$ as follows
		\begin{align*}
			&	\int_{\sum_{i=2}^k \lambda_i^2>2x^2>2} J_k(\lambda_1,\ldots, \lambda_k)d\lambda_1\cdots d\lambda_k  \leq C\int_0^{\infty} (|\lambda_1|+1)^{k-1}\exp(-\lambda_1^2/4) d\lambda_1 \\
			&	\times 	\int_{\sum_{i=2}^k \lambda_i^2>2x^2}  \left(\sum_{i=2}^k \lambda_i^2 \right)^{k(k-1)/4} \exp\left(-\sum_{i=2}^k \lambda_i^2/4 \right) d\lambda_2\cdots d\lambda_k. 
		\end{align*}
		The first integral is bounded.  Using the polar coordinate to the second integral, and by \eqref{r^2} we have the bound 
		\begin{equation}\label{w4}
			I_2\leq C\int_{r>\sqrt{2}x}r^{k(k-1)/2}\exp(-r^2/4)r^{k-2}dr\leq C x^{k(k+1)/2-3}\exp(-x^2/2),
		\end{equation}
		which can be further bounded by $C'x^{k-2}\exp(-x^2/4)$ for $x>1$ by choosing $C'$ large enough.


		Now we estimate $I_1$ for $x>1$. For the case $x^2-(y+1)^2/2>1$, by \eqref{r^2} and \eqref{updensity}, we have  
		\begin{align*} &\P\left(2\lambda_1^2>(2x^2-(y+1)^2)_+, \,\,y^2\leq \sum_{i=2}^k \lambda_i^2\leq (y+1)^2\right)\\
			\leq &C \int_{\lambda_1^2>x^2-(y+1)^2/2>1} |\lambda_1|^{k-1} \exp(-\lambda_1^2/4)d\lambda_1\\
			& \times \int_{\lambda_2>\cdots >\lambda_k:(y+1)^2\geq \sum_{i=2}^k\lambda_i^2\geq y^2}
			\prod_{2\leq i\leq k} (\abs{\lambda_i}+1) \prod_{2\leq i<j\leq k}(\lambda_i-\lambda_j) \exp\big(-
			\sum_{i=2}^k\lambda_i^2/4\big)d\lambda_2\cdots d\lambda_k  \\
			\leq &Cx^{k-2}\exp(-(x^2-(y+1)^2/2)/4)\\&
			\times \int_{\mathbb R^{k-1}: (y+1)^2\geq \sum_{i=2}^k\lambda_i^2\geq y^2}\prod_{2\leq i\leq k} (\abs{\lambda_i}+1)   \prod_{2\leq i<j\leq k}(|\lambda_i|+|\lambda_j|)\exp(-
			\sum_{i=2}^k\lambda_i^2/4)d\lambda_2\cdots d\lambda_k. 
		\end{align*} 
		We denote the last integral as $$A_k(y):= \int_{\mathbb R^{k-1}: (y+1)^2\geq \sum_{i=2}^k\lambda_i^2\geq y^2}\prod_{2\leq i\leq k} (\abs{\lambda_i}+1)    \prod_{2\leq i<j\leq k}(|\lambda_i|+|\lambda_j|)\exp(-
		\sum_{i=2}^k\lambda_i^2/4)d\lambda_2\cdots d\lambda_k. $$
		For the case $x^2-(y+1)^2/2\leq 1$,  by \eqref{stis} and the arguments as above, we simply have\begin{align*} &\P\left(2\lambda_1^2>(2x^2-(y+1)^2)_+, \,\,y^2\leq \sum_{i=2}^k \lambda_i^2\leq (y+1)^2\right)\\ \leq &C \Big(\int_{\mathbb R} (|\lambda_1|+1)^{k-1} \exp(-\lambda_1^2/4)d\lambda_1\Big) A_k(y)\\ \leq &C A_k(y)\\  \leq & C x^{k-2}\exp(-(x^2-(y+1)^2/2)/4)A_k(y)
		\end{align*} for $x>1$. 
		Therefore, in both cases, by the polar coordinate, we further have the   upper bound, 
		\begin{align*} & \P\left(2\lambda_1^2>(2x^2-(y+1)^2)_+, \,\,y^2\leq \sum_{i=2}^k \lambda_i^2\leq (y+1)^2\right)\\ \leq &Cx^{k-2}e^{-x^2/4} \int_{y+1>r>y}(1+r)^{k-1}r^{(k-1)(k-2)/2} r^{k-2}\exp((r+1)^2/8-r^2/4)dr.
		\end{align*}By taking the summation, $I_1$ can be bounded from above by 
		\begin{align*}
			I_1\leq & Cx^{k-2}\exp(-x^2/4) \int_0^{\lfloor \sqrt{2}x\rfloor+1}(1+r)^{k-1}r^{(k-1)(k-2)/2} r^{k-2}\exp((r+1)^2/8-r^2/4)dr\\
			\leq  & Cx^{k-2}\exp(-x^2/4) \int_0^{+\infty} \dots =C' x^{k-2}\exp(-x^2/4).
		\end{align*}
		This will complete the proof of \eqref{w3} by the estimates of $I_1$ and $I_2$.
	\end{proof}
	\subsection{Proof of Lemma \ref{lem:WS1}}
	\begin{proof}
		We first prove \eqref{alpha}. 
		We  claim that \eqref{alpha} is equivalent to the following statement: for any $\delta>0$, $x>1$, there exists a constant $C$ depending on $\delta$, such that
		\begin{align}
			\label{alpha2} &\mathbb{P}(\lambda_1({G}_{\alpha})> x|{G}_{\beta}, {G}_{\gamma})\mathbf{1}_{\{\lambda_1^*(G_{\beta})\leq (1-\delta)x,\lambda_1^*(G_{\gamma})\leq (1-\delta)x\}} \\
			\leq& Cx^{-1}(x/(\lambda_1^*(G_{\beta})+1)+x/(\lambda_1^*(G_{\gamma})+1))^{k(m-k)-1}e^{-(x-\lambda_1^*(G_{\beta}))(x-\lambda_1^*(G_{\gamma}))/2}.\nonumber \end{align} 
		The implication that \eqref{alpha} $\Rightarrow$ \eqref{alpha2} is trivial. We now show that \eqref{alpha2} implies \eqref{alpha}. For any $\delta,\delta'\in (0,1)$, we define $\bar{\delta}=\min\{\delta,\delta'\}$.
		\eqref{alpha2} implies that there exists a constant $C(\bar \delta)$ such that
		\begin{align}
			\label{alpha3} &\mathbb{P}(\lambda_1({G}_{\alpha})> x|{G}_{\beta}, {G}_{\gamma})\mathbf{1}_{\{\lambda_1^*(G_{\beta})\leq (1-\bar \delta)x,\lambda_1^*(G_{\gamma})\leq (1-\bar \delta)x\}} \\
			\leq& Cx^{-1}(x/(\lambda_1^*(G_{\beta})+1)+x/(\lambda_1^*(G_{\gamma})+1))^{k(m-k)-1}e^{-(x-\lambda_1^*(G_{\beta}))(x-\lambda_1^*(G_{\gamma}))/2}.\nonumber \end{align} 
		\eqref{alpha3} implies \eqref{alpha}  since 
		\[
		\Big\{ \lambda_1^*(G_{\beta})\leq (1- \delta)x,\lambda_1^*(G_{\gamma})\leq (1- \delta')x \Big\} 
		\subset
		\Big\{ \lambda_1^*(G_{\beta})\leq (1-\bar \delta)x,\lambda_1^*(G_{\gamma})\leq (1-\bar \delta)x \Big\}. 
		\]
		This completes the proof of the equivalence between \eqref{alpha} and \eqref{alpha2}. We now prove \eqref{alpha2}. 
		Without loss of generality, we may  assume $$\alpha=\{1,\ldots, m\},\gamma=\{1,\ldots, k\},\beta=\{k+1,\ldots, m\}.$$ 
		Since $G_{\beta}$ and $G_{\gamma}$ are both symmetric matrices sampled from GOE (independently), we can find  orthogonal matrices $U$ and $U'$ such that 
		\begin{equation*}
			U G_{\beta}U^t=X, \,\,\,\,
			U' G_{\gamma}{U'}^t=Z,
		\end{equation*}
		where the diagonal matrices 
		\[
		X=\begin{pmatrix}
			\lambda_1(G_{\beta})& & &\\
			& \lambda_2(G_{\beta}) & &\\
			& & \ddots & \\
			& & & \lambda_{\ell}(G_\beta)
		\end{pmatrix},
		Z=\begin{pmatrix}
			\lambda_{1}(G_\gamma)& & &\\
			& \lambda_2(G_{\gamma}) & &\\
			& & \ddots & \\
			& & & \lambda_{k}(G_\gamma) \end{pmatrix},
		\]where  $ \ell:=m-k$.
		
		It follows that
		\[
		\begin{pmatrix}
			U &  \\
			& U'
		\end{pmatrix}
		G_{\alpha}
		\begin{pmatrix}
			U^t &  \\
			& {U'}^t
		\end{pmatrix}
		= \begin{pmatrix}
			X & V \\
			V^t & Z
		\end{pmatrix},
		\]
		where $V$ is an $\ell\times k$ matrix with i.i.d. $N_{\mathbb{R}}(0,1)$ entries.  
		
		Given any $1\times m$ vector $\mathbf v$, we decompose it as 
		$$
		\mathbf v= (\mathbf{p},\mathbf{q}),\,\,\,\mathbf{p}=(p_1,\ldots, p_{\ell}), \,\,\,\mathbf{q}=(q_1,\ldots, q_k), $$
		then we have
		\[
		\mathbf{v}\begin{pmatrix}
			X &  V\\
			V^t & Z
		\end{pmatrix}\mathbf{v}^t
		=\sum_{i=1}^{\ell}\lambda_i(G_{\beta})p_i^2+\sum_{j=1}^{k}\lambda_{j}(G_\gamma)q_j^2+2\sum_{i=1}^{\ell}\sum_{j=1}^k p_iv_{ij}q_j.
		\]
		For simplicity, we use $\P^*$ to denote the conditional probability (conditional on $G_{\beta}$ and $G_{\gamma}$).
		By Rayleigh quotient, for $x>1$,  we have 
		\[
		\begin{split}
			&\P^*(\lambda_1(G_{\alpha})>x)\\
			=&\P^*\left(\exists (\mathbf{p},\mathbf{q})\neq \mathbf{0}, \sum_{i=1}^{\ell}\lambda_i(G_{\beta})p_i^2+\sum_{j=1}^{k}\lambda_{j}(G_\gamma)q_j^2+2\sum_{i=1}^{\ell}\sum_{j=1}^k p_iv_{ij}q_j \geq x\sum_{i=1}^{\ell}p_i^2+x\sum_{j=1}^k q_j^2\right)\\
			=&\P^*\left(\exists (\mathbf{p},\mathbf{q})\neq \mathbf{0}, 2\sum_{i=1}^{\ell}\sum_{j=1}^k p_iv_{ij}q_j \geq \sum_{i=1}^{\ell}(x-\lambda_i(G_{\beta}))p_i^2+\sum_{j=1}^k (x-\lambda_{j}(G_\gamma) q_j^2\right).
		\end{split}
		\]We define the event $$\Omega:=\Big\{\lambda_1^*(G_{\beta})\leq (1-\delta)x,\lambda_1^*(G_{\gamma})\leq (1-\delta)x\Big\}.$$ For the rest of the proof, all the arguments are restricted on the event $\Omega$ and $x>1$.  
		
		On $\Omega$,  by the definition $\lambda^*(A)$ for any symmetric matrix $A$ in \eqref{trace12},  we simply have $x+\sqrt 2(1-\delta)x\geq x-\lambda_i(G_{\beta})\geq x-\lambda_1^*(G_{\beta})\geq \delta x>0$ and $x+\sqrt 2(1-\delta)x\geq x-\lambda_{j}(G_\gamma)\geq x-\lambda_1^*(G_{\gamma})\geq \delta x>0$. 
		
		Therefore,  on $\Omega$ it holds that 
		\begin{equation}\label{zij}
			z_{ij}:=(x-\lambda_i(G_{\beta}))(x-\lambda_{j}(G_\gamma))\in [\delta^2 x^2,  (\sqrt 2-\sqrt2 \delta+1)^2  x^2]. 
		\end{equation} 
		We  denote $$ 
		\bar p_i=\sqrt{x-\lambda_i(G_{\beta})}p_i,\quad 
		\bar q_j=\sqrt{x-\lambda_{j}(G_\gamma)}q_j,\quad \bar v_{ij}=\frac{v_{ij}}{\sqrt{z_{ij}}}.$$
		Using Cauchy-Schwartz inequality we have
		\[
		\begin{split}
			2\sum_{i=1}^{\ell}\sum_{j=1}^k \bar p_i\bar v_{ij} \bar q_j&\leq 2\sqrt{
				\left(\sum_{i=1}^{\ell} \bar p_i^2  \right)
				\left(\sum_{i=1}^{\ell}\left(\sum_{j=1}^k \bar v_{ij}\bar q_j \right)^2\right)
			}\\
			&\leq 2\sqrt{
				\left(\sum_{i=1}^{\ell} \bar p_i^2  \right)
				\left(\sum_{j=1}^k \bar q_j^2 \right)
				\left(\sum_{i=1}^{\ell}\sum_{j=1}^k \bar v_{ij}^2\right)}\\
			&\leq \sqrt{ \sum_{i=1}^{\ell}\sum_{j=1}^k \bar v_{ij}^2}
			\left( \sum_{i=1}^{\ell} \bar p_i^2+\sum_{j=1}^k \bar q_j^2    \right).
		\end{split}
		\]
		Therefore, we further have the probabilistic estimate of $\lambda_1(G_\alpha)$ as
		\[
		\begin{split}
			&\P^*(\lambda_1(G_{\alpha})>x)\\
			=&
			\P^*\left(\exists (\mathbf{\bar p},\mathbf{\bar q})\neq \mathbf{0}, 2\sum_{i=1}^{\ell}\sum_{j=1}^k \bar p_i\bar v_{ij} \bar q_j \geq \sum_{i=1}^{\ell}\bar p_i^2+\sum_{j=1}^k \bar q_j^2\right)\\
			\leq & \P^*\left( \sum_{i=1}^{\ell}\sum_{j=1}^k \bar v_{ij}^2\geq 1 \right).
		\end{split}
		\]
		If $k=\ell=1$,  on the event $\Omega$, by \eqref{r^2} for Gaussian random variable $v_{11}$ we have 
		\begin{equation}
			\begin{split}
				\P^*(\bar v_{11}^2\geq 1)&=\P^*(v_{11}^2>z_{11})\leq \frac{C}{\sqrt {z_{11}}}\exp(-z_{11}/2)\\
				&\leq \frac{C}{x}\exp(-(x-\lambda_1(G_{\beta}))(x-\lambda_{1}(G_\gamma))/2)\\ &\leq \frac{C}{x}\exp(-(x-\lambda_1^*(G_{\beta}))(x-\lambda_1^*(G_{\gamma}))/2),
			\end{split}
		\end{equation}
		where we have used the fact that $z_{11}\geq \delta^2 x^2$. This will imply \eqref{alpha2} in the case $m=2$, $k=\ell=1$. 
		Now we consider the case $k\geq 2$ or $\ell\geq 2$. We define
		\[
		\kappa=\
		\begin{cases}
			\min\left\{\frac{\lambda_1(G_{\beta})-\lambda_2(G_{\beta})}{x-\lambda_2(G_{\beta})},
			\frac{\lambda_{1}(G_\gamma)-\lambda_2(G_{\gamma})}{x-\lambda_2(G_{\gamma})}
			\right\} &\text{if $\ell,k\geq 2$;} \\
			\frac{\lambda_1(G_{\beta})-\lambda_2(G_{\beta})}{x-\lambda_2(G_{\beta})}
			&\text{if $k=1, \ell\geq 2$;}\\
			\frac{\lambda_1(G_{\gamma})-\lambda_2(G_{\gamma})}{x-\lambda_2(G_{\gamma})}
			&\text{if $\ell=1, k\geq 2$.} 
		\end{cases}
		\]
		On the event $\Omega$, one can show that \begin{equation}\label{kap}0 \leq \kappa \leq \frac{(1-\delta)x- \lambda_2(G_{\beta})}{x-\lambda_2(G_{\beta})}\leq  \frac{(1-\delta)x- (-\sqrt 2(1-\delta)x)}{x-(-\sqrt 2(1-\delta)x)}=\frac{(1+\sqrt{2})(1-\delta)}{1+\sqrt 2(1-\delta)}:=c_\delta<1. \end{equation} Note that 
		\[
		\begin{split}
			\sum_{i=1}^{\ell}\sum_{j=1}^k \bar v_{ij}^2=&
			\sum_{i=1}^{\ell}\sum_{j=1}^k  \frac{v_{ij}^2}{(x-\lambda_i(G_{\beta}))(x-\lambda_{j}(G_\gamma))}\\
			= &\frac{1}{(x-\lambda_1(G_{\beta}))(x-\lambda_{1}(G_\gamma))}
			\left(
			v_{11}^2+ \sum_{(i,j)\neq (1,1)}\frac{
				(x-\lambda_1(G_{\beta}))(x-\lambda_{1}(G_\gamma)) 
			}{(x-\lambda_i(G_{\beta}))(x-\lambda_{j}(G_\gamma)) }v_{ij}^2
			\right)\\
			\leq &\frac{1}{(x-\lambda_1(G_{\beta}))(x-\lambda_{1}(G_\gamma))} \left(v_{11}^2+ (1-\kappa)\sum_{(i,j)\neq (1,1)}v_{ij}^2 \right),
		\end{split}
		\]where we used the fact that 
		\begin{equation*}
			\begin{split}
				\frac{
					(x-\lambda_1(G_{\beta}))(x-\lambda_{1}(G_\gamma)) 
				}{(x-\lambda_i(G_{\beta}))(x-\lambda_{j}(G_\gamma)) }&\leq \max \left\{  \frac{
					x-\lambda_1(G_{\beta})}{x-\lambda_2(G_{\beta})},\frac{
					x-\lambda_{1}(G_\gamma)}{x-\lambda_2(G_{\gamma})} \right\}\\
				&=1-\min \left\{ \frac{
					\lambda_1(G_{\beta})-\lambda_2(G_{\beta})}{x-\lambda_2(G_{\beta})},\frac{
					\lambda_{1}(G_\gamma)-\lambda_2(G_{\gamma})}{x-\lambda_2(G_{\gamma})} \right\}
			\end{split}
		\end{equation*}
		$$$$ for $(i,j)\neq (1,1)$. 
		Hence, by \eqref{xi^2} we can conclude  the following bounds on the event $\Omega$ 
		\begin{equation}\label{easybd}
			\begin{split}&\P^*(\lambda_1(G_{\alpha})>x) \\ \leq &  \P^*\left( \sum_{i=1}^{\ell}\sum_{j=1}^k \bar v_{ij}^2\geq 1 \right) \\
				\leq& \P^*\left(v_{11}^2+(1-\kappa)\sum_{(i,j)\neq (1,1)}v_{ij}^2\geq z_{11}\right) \\
				\leq &\P^*\left(\sum_{i=1}^{\ell}\sum_{j=1}^k v_{ij}^2\geq z_{11}\right) \\
				\leq & Cz_{11}^{k\ell/2-1}\exp(-z_{11}/2)\\ 
				\leq &  Cx^{k\ell-2}\exp(-(x-\lambda_1(G_{\beta}))(x-\lambda_{1}(G_\gamma))/2),
			\end{split}
		\end{equation} where in the last step we used \eqref{zij}, and $C$ is a constant depending on $\delta$ and $m$. 
		
		We now consider the case when one of the following two conditions holds. 
		
		Condition 1: $\lambda_1^*(G_{\beta})\leq 1$ or $\lambda_1^*(G_{\gamma})\leq 1$. 
		Under this condition, 
		since $\lambda_1(G_{\beta})\leq \lambda_1^*(G_{\beta})$ and $\lambda_{1}(G_\gamma)\leq \lambda_1^*(G_{\gamma})$, we have 
		$$
		\P^*\left( \sum_{i=1}^{\ell}\sum_{j=1}^k \bar v_{i,j}^2\geq 1 \right) \leq 
		Cx^{k\ell-2}\exp(-(x-\lambda_1^*(G_{\beta}))(x-\lambda_1^*(G_{\gamma}))/2),
		$$
		which  further yields the bound \eqref{alpha2} by the assumption $\lambda_1^*(G_{\beta})\leq 1$ or $\lambda_1^*(G_{\gamma})\leq 1$. 
		
		Condition 2: $
		\max\{\abs{\lambda_i(G_{\beta})},2\leq i\leq \ell\}>\abs{\lambda_1(G_{\beta})}/2
		\mbox{ or }
		\max\{\abs{\lambda_{i}(G_\gamma)},2\leq i\leq k\}>\abs{\lambda_{1}(G_\gamma)}/2.$
		If  condition 1 fails (i.e., $\lambda_1^*(G_{\beta})> 1$ and $\lambda_1^*(G_{\gamma})> 1$) and condition 2 holds (say, $ \max\{\abs{\lambda_i(G_{\beta})},2\leq i\leq \ell\}>\abs{\lambda_1(G_{\beta})}/2$), then we have 
		$$\lambda_1^*(G_{\beta})-\lambda_1(G_{\beta})=\frac{\sum_{i=2}^{\ell}\lambda^2_{i}(G_\beta) }{2(\lambda_1^*(G_{\beta})+\lambda_1(G_{\beta}))}\geq \frac{C(\lambda_1^*(G_{\beta}))^2}{4\lambda_1^*(G_{\beta})}= C'\lambda_1^*(G_{\beta}) .
		$$
		In this case, on $\Omega$ we have
		\begin{equation}\label{gap}
			\begin{split}
				(x-\lambda_1^*(G_{\beta}))(x-\lambda_1^*(G_{\gamma}))&\leq  (x-\lambda_1(G_{\beta}))(x-\lambda_1^*(G_{\gamma}))-(\lambda_1^*(G_{\beta})-\lambda_1(G_{\beta}))\delta x\\
				&\leq  (x-\lambda_1(G_{\beta}))(x-\lambda_{1}(G_\gamma))-C' \lambda_1^*(G_{\beta})\delta x
				\\
				&\leq  (x-\lambda_1(G_{\beta}))(x-\lambda_{1}(G_\gamma))-C'  \delta x .
			\end{split}
		\end{equation}
		Using \eqref{easybd}, we  have 
		\[
		\P^*(\lambda_1(G_{\alpha})>x)\leq Cx^{k\ell-2}\exp(- (x-\lambda^*_1(G_{\beta}))(x-\lambda_1^*(G_{\gamma}))/2) \exp(-C'\delta x).
		\] This can be further bounded from above by 
		$$
		Cx^{-1}\Big(\frac x{1+\lambda^*_1(G_{\beta})}+\frac x{1+\lambda^*_1(G_{\gamma})}\Big)^{k\ell -1}\exp(-(x-\lambda^*_1(G_{\beta}))(x-\lambda_1^*(G_{\gamma}))/2),
		$$
		this is because on $\Omega$ it holds $$\exp(-C'\delta x)\leq  C\Big(\frac 1{1+ (1-\delta)x}\Big )^{k\ell -1}  \leq C\Big(\frac 1{1+\lambda^*_1(G_{\beta})}+\frac 1{1+\lambda^*_1(G_{\gamma})}\Big)^{k\ell -1}$$ for $x>1$ by choosing $C$ large enough.  Note that  the constant $C>0$ only depends on $m$ and $\delta$, and does not depend on (the conditional) $\lambda^*_1(G_{\beta})$ and $\lambda^*_1(G_{\gamma})$.
		
		As a summary,  we have verified \eqref{alpha2} when either of the two conditions is satisfied. 
		Now we assume that both conditions fail so that $\lambda_1^*(G_{\beta})>1$, $ \lambda_1^*(G_{\gamma})>1$, $\max\{\abs{\lambda_i(G_{\beta})},2\leq i\leq \ell\}\leq \abs{\lambda_1(G_{\beta})}/2$ and $\max\{\abs{\lambda_{i}(G_\gamma)},2\leq i\leq k\}\leq \abs{\lambda_{1}(G_\gamma)}/2.$ Then there exists a constant $c>0$ such that
		\[
		\lambda_1(G_{\beta})-\lambda_2(G_{\beta})\geq c\lambda_1^*(G_{\beta})\geq c \mbox{ and }\lambda_{1}(G_\gamma)-\lambda_2(G_{\gamma})\geq c\lambda_1^*(G_{\gamma})\geq c.
		\] 
		Therefore,  by definition of $\kappa$, \eqref{kap} and the assumptions $\lambda_1^*(G_{\beta})>1$, $ \lambda_1^*(G_{\gamma})>1$, we have 
		\begin{equation}\label{kappa1}
			1> c_\delta\geq \kappa \geq c\min\Big\{\frac{\lambda_1^*(G_{\beta})}{x-\lambda_2(G_{\beta})},\frac{\lambda_1^*(G_{\gamma})}{x-\lambda_2(G_{\gamma})}  \Big\}\geq  \frac c2\min\Big\{\frac{\lambda_1^*(G_{\beta})+1}{x-\lambda_2(G_{\beta})},\frac{\lambda_1^*(G_{\gamma})+1}{x-\lambda_2(G_{\gamma})}  \Big\}. 
		\end{equation}
		By the fact that $0<\delta x\leq  x-\lambda_2(G_{\beta}), x-\lambda_2(G_{\gamma}) \leq x+\sqrt 2(1-\delta)x$ on $\Omega$, we further have 
		\begin{equation}\label{kappa2}
			c_\delta\geq  \kappa \geq C_\delta \min\Big\{\frac{\lambda_1^*(G_{\beta})+1}{x},\frac{\lambda_1^*(G_{\gamma})+1}{x}  \Big\}, 
		\end{equation} which implies 
		\begin{equation}\label{1/kappa}
			\frac{1}{\kappa}\leq Cx\left( \frac{1}{\lambda_1^*(G_{\beta})+1}+\frac{1}{\lambda_1^*(G_{\gamma})+1} \right) 
		\end{equation}
		for some constant $C>0$ that depends on $m$ and $\delta$.
		Note that the probability density function $p(t)$ for $v_{11}^2$ is $\exp(-t/2)/\sqrt{2\pi t}$, recall \eqref{easybd},  we have 
		\[
		\begin{split}
			& \P^*\left( \sum_{i=1}^{\ell}\sum_{j=1}^k \bar v_{ij}^2\geq 1 \right)\\
			\leq &\P^*\left(v_{11}^2+(1-\kappa)\sum_{(i,j)\neq (1,1)}v_{ij}^2>z_{11}\right)\\
			= &  \int_{0}^{z_{11}} \frac{\exp(-(z_{11}-t)/2)}{\sqrt{2\pi (z_{11}-t)}}\P^*\left(
			(1-\kappa)\sum_{(i,j)\neq (1,1)}v_{ij}^2>t
			\right) dt+\P^*(v_{11}^2>z_{11})\\
			\leq & C\exp(-z_{11}/2)\left(\frac{1}{\sqrt{z_{11}}}+ \int_0^{z_{11}} \frac{\exp(t/2)}{\sqrt{z_{11}-t}}\P^*\left(
			(1-\kappa)\sum_{(i,j)\neq (1,1)}v_{ij}^2>t
			\right)dt \right),
		\end{split}
		\]
		where we have used 
		$\P^*(v_{11}^2>z_{11})\leq \frac{Ce^{-z_{11}/2}}{\sqrt{z_{11}}}$ by \eqref{xi^2}. 
		By \eqref{xi^2} again, we have 
		\[
		\begin{split}
			\P^*\left( (1-\kappa)\sum_{(i,j)\neq (1,1)}v_{i,j}^2>t\right)&\leq C\left(\frac{t}{1-\kappa}\right)^{(k\ell-1)/2-1}\exp\left(-\frac{t}{2(1-\kappa)}\right)  \\ & \leq Ct^{(k\ell-1)/2-1}\exp\left(-\frac{t}{2(1-\kappa)}\right),
		\end{split}
		\]
		where in the last inequality we used the fact that $1-\kappa\geq 1-c_\delta>0 $ by \eqref{kap}.
		It follows that
		\[
		\begin{split}
			& \P^*\left( \sum_{i=1}^{\ell}\sum_{j=1}^k \bar v_{ij}^2\geq 1 \right)\\
			\leq & C\exp(-z_{11}/2)\left(\frac{1}{\sqrt{z_{11}}}+  \int_0^{z_{11}} \frac{t^{(k\ell-3)/2}}{\sqrt{z_{11}-t}}\exp\left(\frac{t}{2}-\frac{t}{2(1-\kappa)} \right)dt
			\right)\\
			\leq  & C'\exp(-z_{11}/2)\left(\frac{1}{\sqrt{z_{11}}}+ \frac{1}{\sqrt{z_{11}}}\int_0^{z_{11}/2} t^{(k\ell-3)/2}\exp\left(-\frac{\kappa t}{2(1-\kappa)}\right) dt \right.\\
			&\left. +  z_{11}^{(k\ell-3)/2}\exp\left(-\frac{\kappa z_{11}}{4(1-\kappa)}\right)\int_{z_{11}/2}^{z_{11}} \frac{1}{\sqrt{z_{11}-t}}dt  \right)\\
			\leq &  C'\frac{\exp(-z_{11}/2)}{\sqrt{z_{11}}}\Big(1+ \left(\frac{1}{\kappa}\right)^{(k\ell-1)/2} \int_0^{\infty} s^{(k\ell-3)/2}\exp(-s ) ds\\ &+ z_{11}^{(k\ell-1)/2}\exp\Big(-\frac{\kappa z_{11}}{4 }\Big)  \Big).
		\end{split}
		\]
		Note that the integration of $s^{(k\ell-3)/2}\exp(-s)$ converges since $k\ell\geq 2$, and thus the second term can be bounded from above by $c(1 /\kappa)^{(k\ell -1)/2}$. For the third term,  the global maxima of the function $x^{(k\ell-1)/2 }e^{-\kappa x/4}$ is obtained at the point $x=c'/\kappa$, here $c'$ depends on $k, \ell$, and thus the third term can be bounded from above by $C (1/\kappa)^{(k\ell-1)/2}$, where $C$ depends on $\delta$ and $m$. Therefore, on $\Omega$ we have \[
		\begin{split}
			& \P^*\left( \sum_{i=1}^{\ell}\sum_{j=1}^k \bar v_{ij}^2\geq 1 \right)\\
			\leq &  C'\frac{\exp(-z_{11}/2)}{\sqrt{z_{11}}}\Big(1+ c(1 /\kappa)^{(k\ell -1)/2} +C (1/\kappa)^{(k\ell-1)/2}\Big)\\
			\leq &  C'\frac{\exp(-z_{11}/2)}{\sqrt{z_{11}}}\Big((1 /\kappa)^{(k\ell -1)/2}+ c(1 /\kappa)^{(k\ell -1)/2} +C (1/\kappa)^{(k\ell-1)/2}\Big)\,\,\, [\mbox{since}\,\,\, 1<1/\kappa]\\ =& C\frac{\exp(-z_{11}/2)}{\sqrt{z_{11}}}(1 /\kappa)^{(k\ell -1)/2}\\ \leq & C\frac{\exp(-z_{11}/2)}{\sqrt{z_{11}}}(1 /\kappa)^{k\ell -1}\\
			\leq & Cx^{-1}  \left(\frac{x}{\lambda_1^*(G_{\beta})+1} +\frac{x}{\lambda_1^*(G_{\gamma})+1} \right)^{k\ell-1} \exp(-(x-\lambda_1^*(G_{\beta}))(x-\lambda_1^*(G_{\gamma}))/2) )\,\,\, [\mbox{by}\,\,\eqref{zij}, \eqref{1/kappa}].
		\end{split}
		\]This proves \eqref{alpha2}, and thus  we complete the proof of \eqref{alpha}.


		We now prove \eqref{alphabeta} using \eqref{w3} and \eqref{alpha}. 
		By convexity of the function $s\to s^{k(m-k)-1}$ for $s>0$, we have 
		\begin{equation}\label{i3i4}
			\begin{split}
				&x^{-1}(x/(\lambda_1^*(G_{\beta})+1)+x/(\lambda_1^*(G_{\gamma})+1))^{k(m-k)-1}e^{-(x-\lambda_1^*(G_{\beta}))(x-\lambda_1^*(G_{\gamma}))/2}\\
				\leq &Cx^{-1}(x/(\lambda_1^*(G_{\beta})+1))^{k(m-k)-1} 
				e^{-(x-\lambda_1^*(G_{\beta}))(x-\lambda_1^*(G_{\gamma}))/2}\\
				&+Cx^{-1}(x/(\lambda_1^*(G_{\gamma})+1))^{k(m-k)-1} 
				e^{-(x-\lambda_1^*(G_{\beta}))(x-\lambda_1^*(G_{\gamma}))/2}\\
				:=&I_3+I_4.
			\end{split}
		\end{equation} Let $\bar \E$ be the conditional expectation with respect to $G_{\gamma}$ and $\ell:=m-k$ as before, and we define
		\[
		f(t)=x^{-1}(x/(t+1))^{k\ell -1}\exp(-(x-t)(x-\lambda_1^*(G_{\gamma}))/2)
		\] and 
		$$
		h(t)=x^{-1}(x/(\lambda_1^*(G_{\gamma})+1))^{k\ell -1}\exp(-(x-t)(x-\lambda_1^*(G_{\gamma}))/2).
		$$
		We define the event 
		$$\Omega' =\{\lambda_1^*(G_{\gamma})\leq (1-\delta')x\}.$$ 
		Using \eqref{alpha} and \eqref{i3i4}, we have 
		\begin{equation}\label{esalpha}
			\begin{split}
				&\mathbb{P}(\lambda_1({G}_{\alpha})> x,\ \lambda_1^*({G}_{\beta})\leq (1-\delta)x| {G}_{\gamma})\mathbf{1}_{\{\lambda_1^*(G_{\gamma})\leq (1-\delta')x\}}\\ 
				=& \E[\mathbb{E}\left(\mathbf 1_{\{\lambda_1({G}_{\alpha})> x,\  {\lambda^*_1(G_{\beta})<(1-\delta) x}\}}|G_\beta, G_\gamma\right)|G_\gamma]\mathbf{1}_{\{\lambda_1^*(G_{\gamma})\leq (1-\delta')x\}}
				\\=& \E[\mathbb{E}\left(\mathbf 1_{\{\lambda_1({G}_{\alpha})> x\}}|G_\beta, G_\gamma\right)\mathbf 1_{\{\lambda^*_1(G_{\beta})<(1-\delta)x\}}\mathbf{1}_{\{\lambda_1^*(G_{\gamma})\leq (1-\delta')x\}}|G_\gamma]\mathbf{1}_{\{\lambda_1^*(G_{\gamma})\leq (1-\delta')x\}}
				\\ 
				\leq &C\bar \E\Big(x^{-1}(x/(\lambda_1^*(G_{\beta})+1)+x/(\lambda_1^*(G_{\gamma})+1))^{k\ell-1}\\ &\times\exp\big(-(x-\lambda_1^*(G_{\beta}))(x-\lambda_1^*(G_{\gamma}))/2\big)
				\1_{\{\lambda_1^*(G_{\beta})\leq (1-\delta)x\}}  \Big)\mathbf{1}_{_{\Omega'}}\\
				\leq &C \bar \E(I_3     \1_{\{\lambda_1^*(G_{\beta})\leq (1-\delta)x\}})\mathbf{1}_{_{\Omega'}}+C\bar \E(I_4     \1_{\{\lambda_1^*(G_{\beta})\leq (1-\delta)x\}})\mathbf{1}_{{\Omega'}}\\
				\leq &C\Big(\bar \E(f(\lambda_1^*(G_{\beta}))\1_{\{\lambda_1^*(G_{\beta})<(1-\delta)x\}})\mathbf{1}_{\Omega'}+ \bar \E(h(\lambda_1^*(G_{\beta}))\1_{\{\lambda_1^*(G_{\beta})<(1-\delta)x\}})\mathbf{1}_{\Omega'}\Big).
			\end{split}
		\end{equation}
		Note that $G_{\beta}$ and $G_{\gamma}$ are independent, therefore, 
		$\lambda_1^*(G_{\beta})$ and $\lambda_1^*(G_{\gamma})$ are independent as well.
		Hence for any $t>0$, $\bar \P(\lambda_1^*(G_{\beta})>t)$ is equal to 
		unconditional probability $\P(\lambda_1^*(G_{\beta})>t)$. 
		Then using integration by parts, we have 
		\begin{equation}\label{parts}
			\begin{split}
				\bar  \E(f(\lambda_1^*(G_{\beta}))\1_{\{\lambda_1^*(G_{\beta})<(1-\delta)x\}})=&\int_0^{(1-\delta)x} f'(t)\P(\lambda_1^*(G_{\beta})>t)dt\\
				&-f((1-\delta)x)\P(\lambda_1^*(G_{\beta})>(1-\delta)x) +f(0).
			\end{split}
		\end{equation}
		On the event $\Omega' =\{\lambda_1^*(G_{\gamma})\leq (1-\delta')x\},$ one has
		$$
		f'(t)\leq \frac{1}{2}   (x/(t+1))^{k\ell -1}\exp(-(x-t)(x-\lambda_1^*(G_{\gamma}))/2).
		$$ 
		Therefore,    we have
		\begin{equation}\label{f'01}
			\begin{split}
				&\int_0^{1}f'(t)\P(\lambda_1^*(G_{\beta})>t)dt\leq \frac{1}{2}x^{k\ell-1}
				\exp(-(x-1)(x-\lambda_1^*(G_{\gamma}))/2)\\
				&=\frac{1}{2}x^{k\ell-1}e^{1/4} \exp((\lambda_1^*(G_{\gamma})^2-x^2 )/4) 
				\exp\left(-\frac{((x-\lambda_1^*(G_{\gamma}))-1)^2}{4}\right).
			\end{split} 
		\end{equation}
		On the event $\Omega'$, we have $x-\lambda_1^*(G_{\gamma})\geq \delta' x$. Therefore, for $C$ large enough, one has
		\begin{equation}\label{x-lambda^*}
			\exp\left(-\frac{((x-\lambda_1^*(G_{\gamma}))-1)^2}{4}\right) \leq C \exp(-\delta'^2 x^2/8).
		\end{equation}
		It further follows from \eqref{f'01} and \eqref{x-lambda^*} that for $x>1$, there exists $C$ large enough such that
		\begin{equation} \label{ddds}
			\int_0^{1}f'(t)\P(\lambda_1^*(G_{\beta})>t)dt\leq C 
			x^{\ell-2}\exp(((\lambda_1^*(G_{\gamma}))^2-x^2 )/4). 
		\end{equation}
		If $(1-\delta)x\leq 1$, then we trivially have 
		\begin{equation}\label{1-deltax<1}
			\int_0^{(1-\delta)x} f'(t)\P(\lambda_1^*(G_{\beta})>t)dt\leq \int_0^{1}\dots\leq C 
			x^{\ell-2}\exp(((\lambda_1^*(G_{\gamma}))^2-x^2 )/4). 
		\end{equation}
		Now we consider the case when $(1-\delta)x\geq 1$. 
		Recall \eqref{w3} in Lemma \ref{lem:WS}, we have 
		\begin{equation}\label{pbetat}
			\P(\lambda_1^*(G_{\beta})>t) \leq C t^{\ell-2}\exp(-t^2/4),\,\,\,t\geq 1,
		\end{equation}
		this implies that
		\begin{equation}\label{f'p}
			\begin{split}
				& \int_1^{(1-\delta)x}f'(t)\P(\lambda_1^*(G_{\beta})>t)dt\leq
				C x^{\ell-2}
				\exp(((\lambda_1^*(G_{\gamma}))^2-x^2 )/4)\\
				&\times \int_1^{(1-\delta)x}     (x/(t+1))^{k\ell -1}\exp(-{(t-(x-\lambda_1^*(G_{\gamma})))^2}/4)dt,
			\end{split}
		\end{equation}
		where we have used the identity \begin{equation*}
			\begin{split}
				& -(x-t)(x-\lambda_1^*(G_{\gamma}))/2-t^2/4\\
				=&-\frac{(t-(x-\lambda_1^*(G_{\gamma})))^2}{4} +\frac{(x-\lambda_1^*(G_{\gamma}))^2}{4}-\frac{x(x-\lambda_1^*(G_{\gamma}))}{2}\\
				=&-\frac{(t-(x-\lambda_1^*(G_{\gamma})))^2}{4}+\frac{(\lambda_1^*(G_{\gamma}))^2-x^2 }{4}.
			\end{split}
		\end{equation*}
		On $\Omega'$,  we have, 
		\begin{equation}\label{ctlint}
			\begin{split}
				&\int_1^{(1-\delta)x}     (x/(t+1))^{k\ell -1}\exp(-{(t-(x-\lambda_1^*(G_{\gamma})))^2}/4)dt\\
				\leq & \int_0^{\min\{\delta' x/2, (1-\delta)x\}}x^{k\ell-1}\exp(-(\delta')^2 x^2/16)dt\\
				& + \int_{\min\{\delta' x/2, (1-\delta)x\}}^{(1-\delta)x}\Big (\frac{2}{\delta'}+\frac{1}{1-\delta} \Big)^{k\ell-1}
				\exp(-{(t-(x-\lambda_1^*(G_{\gamma})))^2}/4)dt\\
				\leq & C x^{k\ell}\exp(-(\delta')^2 x^2/16)+C \leq C'.
			\end{split}
		\end{equation} 
		Combining \eqref{ddds}, \eqref{f'p} and \eqref{ctlint}, for $(1-\delta)x>1$, we have 
		\begin{equation}\label{ctlint2}
			\int_0^{(1-\delta)x}f'(t)\P(\lambda_1^*(G_{\beta})>t)dt=\int_{0}^1\dots+\int_{1}^{(1-\delta)x}\dots \leq
			C x^{\ell-2}
			\exp(((\lambda_1^*(G_{\gamma}))^2-x^2 )/4).
		\end{equation}
		By \eqref{1-deltax<1}, the above estimate is actually true for all $x>1$ on $\Omega'$. 
		We also have
		\begin{equation}\label{fp}
			f(0)=x^{k\ell-2} \exp(-x(x-\lambda_1^*(G_{\gamma}))/2) \leq C x^{\ell-2} \exp(((\lambda_1^*(G_{\gamma}))^2-x^2 )/4)
		\end{equation}
		for $C$ large enough, this is because on $\Omega'$ we have 
		$$  -x(x-\lambda_1^*(G_{\gamma}))/2-((\lambda_1^*(G_{\gamma}))^2-x^2 )/4=-(x-\lambda_1^*(G_{\gamma}))^2/4\leq -\delta'^2x^2.
		$$
		Combining \eqref{parts}, \eqref{ctlint2} and \eqref{fp}, on $\Omega'$ we get
		\begin{equation}\label{ctli3}
			\bar \E(f(\lambda_1^*(G_{\beta}))\1_{\{\lambda_1^*(G_{\beta})<(1-\delta)x\}})\leq C' x^{\ell-2}\exp(((\lambda_1^*(G_{\gamma}))^2-x^2 )/4). 
		\end{equation} 
		On $\Omega'$,  for $x>1$, it holds that $ x /(\lambda^*_{1}(G_\gamma)+1)\geq x/((1-\delta')x+1)\geq 1/(2-\delta')$, therefore, 
		we further have the upper bound,
		\begin{equation}\label{ctli6}
			\bar \E(f(\lambda_1^*(G_{\beta}))\1_{\{\lambda_1^*(G_{\beta})<(1-\delta)x\}})\leq C' x^{\ell-2} (x /(\lambda^*_{1}(G_ \gamma)+1))^{k\ell -1}\exp(((\lambda_1^*(G_{\gamma}))^2-x^2 )/4). 
		\end{equation} 
		
		We can similarly control the conditional expectation of $h(\lambda_1^*(G_{\beta}))\1_{\{\lambda_1^*(G_{\beta})\leq (1-\delta)x\}}$.   Analogously to \eqref{parts}, we have 
		\begin{equation}
			\begin{split}
				\bar  \E(h(\lambda_1^*(G_{\beta}))\1_{\{\lambda_1^*(G_{\beta})<(1-\delta)x\}})=&\int_0^{(1-\delta)x} h'(t)\P(\lambda_1^*(G_{\beta})>t)dt\\
				&-h((1-\delta)x)\P(\lambda_1^*(G_{\beta})>(1-\delta)x) +h(0).
			\end{split}
		\end{equation}
		On the event $\Omega' =\{\lambda_1^*(G_{\gamma})\leq (1-\delta')x\},$ one has
		$$
		h'(t)\leq \frac{1}{2}   (x/(\lambda_1^*(G_{\gamma})+1))^{k\ell -1}\exp(-(x-t)(x-\lambda_1^*(G_{\gamma}))/2).
		$$ We basically repeat the proofs of \eqref{ctlint2} and \eqref{fp} to get
		\begin{equation}\label{ctli4}
			\begin{split}
				&\bar \E(h(\lambda_1^*(G_{\beta}))\1_{\{\lambda_1^*(G_{\beta})<(1-\delta)x\}})
				\leq C x^{\ell-2} (x/(\lambda_1^*(G_{\gamma})+1))^{k\ell-1}\exp(((\lambda_1^*(G_{\gamma}))^2-x^2 )/4).
			\end{split}
		\end{equation}
		Then \eqref{alphabeta}  follows from \eqref{esalpha}, \eqref{ctli6} and \eqref{ctli4} (recall  that $\ell=m-k$). 
	\end{proof}
	
	\subsection{Proof of Lemma  \ref{t_nlimit}}\label{lemma44}
	\begin{proof}
		By  the definition of $A_\alpha$ in \eqref{defAS} with $\alpha=\{1, 2,..., m\}$, we have 
		\begin{equation}\label{a_alpha}
			\begin{split}
				&\mathbb P( A_{\{1, 2,..., m \}})\\=&\mathbb P(\lambda_1({G}_\alpha)>y_{m};\ \ \lambda_1^*({G}_{\beta})\leq y_{k},\ \forall\ 1\leq k<m,\ \beta\subset \alpha,\ |\beta|=k)\\
				=&\P(\lambda_1({G}_\alpha)>y_{m})  \P(\lambda_1^*({G}_{\beta})\leq y_{k},\ \forall\ 1\leq k<m,\ \beta\subset \alpha,\ |\beta|=k|\lambda_1({G}_\alpha)>y_{m}).
			\end{split}
		\end{equation}
		Lemma \ref{t_nlimit}  follows from the following two limits,
		\begin{equation}\label{lam>ym}
			\begin{split}&\lim_{n\to+\infty}\mathbb{P}(\lambda_1(G_{\alpha})>y_m)/  [m2^{-(1+m)/2}/\Gamma(1+m/2)\cdot y_m^{m-2}e^{-y_m^2/4} ]=1
			\end{split}
		\end{equation}
		and \begin{equation}\label{s1}
			\lim_{n\to+\infty}\P(\lambda_1^*({G}_{\beta})\leq y_{k},\ \forall\ 1\leq k<m,\ \beta\subset \alpha,\ |\beta|=k|\lambda_1({G}_\alpha)>y_{m})
			=K_{m}, 
		\end{equation} 
		where $K_m$ is the constant defined in Theorem \ref{main}.  Indeed, \eqref{a_alpha}, \eqref{lam>ym}, \eqref{s1} and the definition of $y_m$ imply  that
		\begin{align*} 
			\lim_{n\to+\infty}t_n&=\lim_{n\to+\infty}{n\choose m}\mathbb{P}(A_{\{1,\cdots,m\}}) \\ &=\lim_{n\to+\infty} n^m/m! \Big(m2^{-(1+m)/2}/\Gamma(1+m/2)\cdot y_m^{m-2}e^{-y_m^2/4}\Big) K_{m}\\&=\lim_{n\to+\infty}\frac{y_m^{m-2}(\ln n)^{-(m-2)/2}K_{m}}{(m-1)!2^{\frac{m+1}{2}}\Gamma(1+m/2)}e^{-y/4}\\&=\frac{(2m)^{(m-2)/2}K_{m}}{(m-1)!2^{3/2}\Gamma(1+m/2)}e^{-y/4},
		\end{align*}
		as desired.

		We first  prove \eqref{lam>ym}.
		By  formula \eqref{jointgoe} for $G_\alpha$ where $|\alpha|=m$,  we have
		\begin{equation}\label{lambda1s1}
			\begin{split}
				&\mathbb{P}(\lambda_1(G_{\alpha})>y_m)\\=&\frac{1}{Z_m}\int_{\lambda_1>\cdots>\lambda_m;\lambda_1>y_m}
				e^{-\sum\limits_{i=1}^m\lambda_i^2/4}\prod_{1\leq i<j\leq m}|\lambda_i-\lambda_j|d\lambda_1\cdots d\lambda_m\\=&\frac{1}{Z_m}\int_{\lambda_1>y_m}
				e^{-\lambda_1^2/4}g(\lambda_1)d\lambda_1,
			\end{split}
		\end{equation}
		where we denote 
		\begin{align*}
			g(\lambda_1)&=\int_{\lambda_1>\cdots>\lambda_m}
			e^{-\sum\limits_{i=2}^m\lambda_i^2/4}\prod_{1\leq i<j\leq m}|\lambda_i-\lambda_j|d\lambda_2\cdots d\lambda_m.\\
		\end{align*}
		We claim that
		\begin{equation}\label{gasymp}
			\begin{split}
				\lim_{\lambda_1\to+\infty }g(\lambda_1)/(\lambda_1^{m-1}Z_{m-1}) =1. 
			\end{split}
		\end{equation}
		Dividing into two cases $\lambda_m>-\sqrt{\lambda_1}$ and $\lambda_m<-\sqrt{\lambda_1}$, $g(\lambda_1)$ is bounded from above by \begin{equation}
			\begin{split}
				&\int_{\lambda_2>\cdots>\lambda_m>-\sqrt{\lambda_1}} (\lambda_1+\sqrt{\lambda_1})^{m-1} \prod_{2\leq i<j\leq m}(\lambda_i-\lambda_j)e^{-\sum_{i=2}^m \lambda_i^2/4}d\lambda_2\cdots d\lambda_m\\
				& +\int_{\lambda_2>\cdots>\lambda_m,\lambda_m<-\sqrt{\lambda_1}}  \prod_{1\leq i<j\leq m}(\lambda_i-\lambda_j)e^{-\sum_{i=2}^m \lambda_i^2/4}d\lambda_2\cdots d\lambda_m\\
				&:=I_5+I_6.
			\end{split}
		\end{equation}
		Note that $I_5$ is further bounded from above by $(\lambda_1+\sqrt{\lambda_1})^{m-1}Z_{m-1}$.
		For $I_6$, we note that 
		$$
		\prod_{1\leq i<j\leq m}(\lambda_i-\lambda_j)\leq
		\prod_{1\leq i<j\leq m}(\abs{\lambda_i}+1) (\abs{\lambda_j}+1)\leq \prod_{i=1}^m (\abs{\lambda_i}+1)^{m-1}.
		$$
		Therefore, we can  bound $I_6$ from above by 
		\begin{equation}\label{43}
			\begin{split}
				&  (\abs{\lambda_1}+1)^{m-1} \int_{\lambda_m<-\sqrt{\lambda_1}}
				(\abs{\lambda_m}+1)^{m-1}\exp(-\lambda_m^2/4)d\lambda_m\\
				\times & \left( \int_{\mathbb{R}}
				(\abs{\lambda_2}+1)^{m-1}\exp(-\lambda_2^2/4)d\lambda_2 \right)^{m-2}.
			\end{split}
		\end{equation}
		For $\lambda_1 $ large enough, \eqref{43} can be further bounded from above by (using \eqref{r^2})
		\[
		C\lambda_1^{m-1} \sqrt{\lambda_1}^{m-2}\exp(-\lambda_1/4)\leq C'\lambda_1^{m-2}.
		\]
		Combining the estimates for $I_5$ and $I_6$ we get 
		\[
		g(\lambda_1)\leq Z_{m-1}(\lambda_1+\sqrt{\lambda_1})^{m-1}+C\lambda_1^{m-2}
		\]
		for $\lambda_1$ large enough.
		It follows that
		\[
		\limsup_{\lambda_1\to+\infty}\frac{g(\lambda_1)}{\lambda_1^{m-1}}\leq Z_{m-1}.
		\]
		For the  lower bound, for all $\lambda_1>1$ we have 
		\[
		g(\lambda_1)\geq (\lambda_1-\sqrt{\lambda_1})^{m-1}
		\int_{\sqrt{\lambda_1}\geq \lambda_2>\cdots>\lambda_m}
		e^{-\sum\limits_{i=2}^m\lambda_i^2/4}\prod_{2\leq i<j\leq m}|\lambda_i-\lambda_j|d\lambda_2\cdots d\lambda_m.
		\]
		It follows that 
		\[
		\liminf_{\lambda_1\to+\infty}\frac{g(\lambda_1)}{\lambda_1^{m-1}}\geq
		\liminf_{\lambda_1\to+\infty}  \int_{\sqrt{\lambda_1}\geq \lambda_2>\cdots>\lambda_m}
		e^{-\sum\limits_{i=2}^m\lambda_i^2/4}\prod_{2\leq i<j\leq m}|\lambda_i-\lambda_j|d\lambda_2\cdots d\lambda_m=Z_{m-1}.
		\]
		Then \eqref{gasymp} follows from the upper and lower bounds. 
		
		Therefore,  as $n\to+\infty$, i.e., $y_m\to+\infty$, we have 
		\begin{equation}\label{star}
			\int_{\lambda_1>y_m}
			e^{-\lambda_1^2/4}g(\lambda_1)d\lambda_1 \sim  Z_{m-1}
			\int_{\lambda_1>y_m} e^{-\lambda_1^2/4} \lambda_1^{m-1}
			d\lambda_1.
		\end{equation}
		Standard results for the upper incomplete Gamma function imply
		\begin{equation}\label{lambda1s2}
			\begin{split}
				\int_{\lambda_1>y_m} e^{-\lambda_1^2/4} \lambda_1^{m-1}
				d\lambda_1 
				&=\int_{s\geq y_m^2/4} e^{-s} 2^{m-1} s^{(m-2)/2} ds \\
				&\sim 2^{m-1} e^{-y_m^2/4} (y_m^2/4)^{(m-2)/2} \\
				&=2e^{-y_m^2/4}y_m^{m-2}.
			\end{split}
		\end{equation}
		Combining \eqref{lambda1s1}, \eqref{star} and \eqref{lambda1s2} we get
		\begin{equation}\label{lam>ym2}
			\begin{split}&\mathbb{P}(\lambda_1(G_{\alpha})>y_m)\\  \sim& \frac{ Z_{m-1}}{Z_m}\int_{\lambda_1>y_m}
				e^{-\lambda_1^2/4}\lambda_1^{m-1}d\lambda_1\\ \sim&( Z_{m-1}/{Z_m})2y_m^{m-2}e^{-y_m^2/4}\\=&m(2\pi)^{-1/2}\Gamma(3/2)/\Gamma(1+m/2)\cdot2^{1-m/2}y_m^{m-2}e^{-y_m^2/4}
				\\=&m2^{-(1+m)/2}/\Gamma(1+m/2)\cdot y_m^{m-2}e^{-y_m^2/4}.
			\end{split}
		\end{equation}
		This completes the proof of \eqref{lam>ym}.

		
		
		Now we   prove \eqref{s1}.  We define the following two auxiliary  events:  
		$$B_{\alpha}:=
		\Big\{y_m<\lambda_1({G}_\alpha)<y_{m}+1,\lambda_2(G_{\alpha})<\sqrt{\log \log n},\lambda_m(G_{\alpha})>-\sqrt{\log \log n}\Big \} 
		$$ and $$
		D_{\alpha}: =\Big\{y_m<\lambda_1(G_{\alpha})<y_m+1\Big\}.
		$$
		Then it holds that $$B_{\alpha}\subset D_{\alpha}\subset \Big\{\lambda_1(G_{\alpha})>y_m\Big\}.$$
		We first show that 
		\begin{equation}\label{s2}
			\lim_{n\to+\infty}\P(B_{\alpha}|\lambda_1(G_{\alpha})>y_m)=1.
		\end{equation} 
		Note that \eqref{s2} follows if we can prove the following two limits, \begin{equation}\label{dalpha}
			\lim_{n\to+\infty} \P(D_{\alpha}|\lambda_1(G_{\alpha})>y_m)=1
		\end{equation}
		and   
		\begin{equation}\label{dtob}
			\lim_{n\to+\infty} \P(B_{\alpha}|D_{\alpha})=1.
		\end{equation}
		By \eqref{lam>ym}, we find that 
		\begin{equation*}
			\lim_{n\to+\infty} \frac{\P(\lambda_1(G_{\alpha})>y_m+1 ) }{\P(\lambda_1(G_{\alpha})>y_m )} =0, 
		\end{equation*}
		which implies \eqref{dalpha}. 
		We now prove \eqref{dtob} by finding 
		a sequence of numbers $\epsilon_n\to 0$ such that  for $n$ large enough  it holds that \begin{equation}\label{s3}
			\frac{\P(\lambda_2(G_{\alpha})>\sqrt{\log \log n} \mbox{ or }\lambda_m(G_{\alpha})<-\sqrt{\log \log n} |D_{\alpha}) }{\P(\lambda_2(G_{\alpha})<1,\lambda_m(G_{\alpha})>0 |D_{\alpha})}\leq \epsilon_n .
		\end{equation}
		Indeed, \eqref{s3} implies that 
		\begin{equation}
			\P(\lambda_2(G_{\alpha})>\sqrt{\log \log n} \mbox{ or }\lambda_m(G_{\alpha})<-\sqrt{\log \log n}|D_{\alpha})\leq \frac{\epsilon_n}{\epsilon_n+1},
		\end{equation}
		which converges to 0 as $n\to+\infty$, and thus this yields \eqref{dtob}. To prove \eqref{s3}, we divide the set 
		$\Big \{\lambda_2(G_{\alpha})>\sqrt{\log \log n} \mbox{ or }\lambda_m(G_{\alpha})<-\sqrt{\log \log n}\Big\}$ as the union of  two disjoint subsets $S_1\cup S_2$ as follows.  On the subset $S_1:=\Big\{\lambda_m(G_{\alpha})\leq -\sqrt{\log \log n}\Big\}$, if we  condition on $D_\alpha$,   we have \begin{equation}
			\begin{split}
				&\prod_{1\leq i< j\leq m}(\lambda_i-\lambda_j) \\ \leq & C (\lambda_1+\abs{\lambda_m})^{m-1}(\abs{\lambda_2}+\abs{\lambda_m})^{(m-1)(m-2)/2}\\
				\leq & C \lambda_1^{m-1}\abs{\lambda_m}^{m-1}
				(\abs{\lambda_2}+1)^{(m-1)(m-2)/2}
				\abs{\lambda_m}^{(m-1)(m-2)/2}\\
				= &C \lambda_1^{m-1}\abs{\lambda_m}^{(m-1)m/2}
				(\abs{\lambda_2}+1)^{(m-1)(m-2)/2}.
			\end{split}
		\end{equation}
		On the subset $S_2:=\Big\{\lambda_2(G_{\alpha})>\sqrt{\log \log n}, \lambda_m(G_{\alpha})>-\sqrt{\log \log n}\Big\}$, if we  condition on $D_\alpha$ where $\lambda_1(G_\alpha)\sim 2\sqrt{m\log n}$, then  we easily have the upper bound
		$$
		\prod_{1\leq i< j\leq m}(\lambda_i-\lambda_j)\leq 
		C\lambda_1^{m-1}\lambda_2^{(m-1)(m-2)/2}.
		$$
		It follows that the
		left hand side of \eqref{s3} can be  bounded from above  by the summation of 
		$$
		\frac{C[\int_{y_m}^{y_m+1}\lambda_1^{m-1}e^{-\lambda_1^2/4}d\lambda_1]\int_{\lambda_m<-\sqrt{\log \log n}} \abs{\lambda_m}^{(m-1)m/2}
			(\abs{\lambda_2}+1)^{(m-1)(m-2)/2} e^{-\sum_{i=2}^m\lambda_i^2/4}d\lambda_2\cdots d\lambda_m }{[\int_{y_m}^{y_m+1}(\lambda_1-1)^{m-1}e^{-\lambda_1^2/4}d\lambda_1]\int_{1>\lambda_2>\cdots >\lambda_m>0} 
			\prod_{2\leq i<j\leq m}(\lambda_i-\lambda_j) 
			e^{-\sum_{i=2}^m\lambda_i^2/4}
			d\lambda_2\cdots d\lambda_m}
		$$
		and 
		$$
		\frac{ C[\int_{y_m}^{y_m+1}\lambda_1^{m-1}e^{-\lambda_1^2/4}d\lambda_1]\int_{\lambda_2>\sqrt{\log \log n}, \lambda_m>-\sqrt{\log \log n}} \lambda_2^{(m-2)(m-1)/2} e^{-\sum_{i=2}^m\lambda_i^2/4}d\lambda_2\cdots d\lambda_m }{[\int_{y_m}^{y_m+1}(\lambda_1-1)^{m-1}e^{-\lambda_1^2/4}d\lambda_1]\int_{1>\lambda_2>\cdots >\lambda_m>0} 
			\prod_{2\leq i<j\leq m}(\lambda_i-\lambda_j) 
			e^{-\sum_{i=2}^m\lambda_i^2/4} d\lambda_2\cdots d\lambda_m}.
		$$
		The above summation can be further bounded from above  by
		\begin{equation*}
			\begin{split}
				&C\Big( \int_{\lambda_m<-\sqrt{\log \log n}} \abs{\lambda_m}^{(m-1)m/2}\exp(-\lambda_m^2/4)d\lambda_m\\
				&+  \int_{\lambda_2>\sqrt{\log \log n}} \lambda_2^{(m-1)(m-2)/2}\exp(-\lambda_2^2/4)d\lambda_2 \Big)\\
				&:=\epsilon_n
			\end{split} 
		\end{equation*}as $n$ large enough. 
		Clearly it holds that $\epsilon_n\to 0$ since $\sqrt{\log \log n}\to+\infty$.  This completes the proof of 
		\eqref{s3} and thus the proof  of \eqref{s2}. 
		
		Now we are ready to prove \eqref{s1}.
		We define the event $$H_{\alpha}=\Big\{\lambda_1^*({G}_{\beta})\leq y_{k},\ \forall\ 1\leq k<m,\ \beta\subset \alpha,\ |\beta|=k\Big\}.$$
		Using \eqref{s2}, we see that \eqref{s1} is equivalent to  
		\begin{equation}\label{s4}
			\lim_{n\to+\infty}\P(H_\alpha |B_{\alpha}) =K_{m}.
		\end{equation}
		To prove \eqref{s4} we need to use the fact that\begin{align*} & {G}_{\alpha}\eqd U^T\diag(\lambda_1,\cdots,\lambda_m)U\,\,
		\end{align*}where $|\alpha|=m$ and $$U\eqd \mathcal{U}(O(m))$$ is sampled from the uniform measure  on  the orthogonal group $O(m)$ which is independent of $(\lambda_1,\cdots,\lambda_m)$ and $$(\lambda_1,\cdots,\lambda_m)\eqd(\lambda_1({G}_{\alpha}),\cdots,\lambda_m({G}_{\alpha})). $$ 
		Given any \begin{align*} &X:= U^T\diag(\lambda_1,\cdots,\lambda_m)U,\ U\in O(m),\ \lambda_1\geq\cdots\geq\lambda_m,
		\end{align*}by definition it holds \begin{align*} &\lambda_1(X)= \lambda_1,\ \lambda_1^*(X)^2=\lambda_1^2+\sum_{k=2}^m\lambda_k^2/2.
		\end{align*}
		Let $\lambda_1^-(X):=\sqrt{\sum_{k=2}^m\lambda_k^2}=\sqrt{|X|^2-\lambda_1^2(X)},\ X_1:=U^T\diag(\lambda_1,0,\cdots,0)U,\ X_2:=X-X_1 ,$ then we have $|X_2| (=\sqrt{\mathrm{ Tr(X_2^2)}} )=\lambda_1^-(X) $. 
		For $ \beta \subset \alpha=\{1,\cdots,m\}$ we have 
		\begin{equation}\label{normineq}
			|\lambda_1^*(X_{\beta})-\lambda_1^*((X_1)_{\beta})|\leq |(X_2)_{\beta}|\leq |X_2|=\lambda_1^-(X).
		\end{equation} 
		Here,   by the Lipschitz continuity of eigenvalues for symmetric matrices (see Corollary A.6 in \cite{AGZ}),  we can further derive the fact that 
		$$|\lambda_1^*(A)-\lambda_1^*(B)|\leq |A-B|.$$
		Since $(X_1)_{\beta}$ is a matrix of rank at most $1$, we have 
		\begin{equation}\label{lambda1*}
			\lambda_1^*((X_1)_{\beta})=(\Tr((X_1)_{\beta}^2))^{1/2}= \abs{\lambda_1}\sum_{k\in\beta}u_k^2,
		\end{equation}
		here $\mathbf u=(u_1,\cdots,u_m)$ is the first row of $U$. Note that $\mathbf u$ has the uniform distribution on the unit sphere $S^{m-1}$.
		By \eqref{normineq} and \eqref{lambda1*}, we have 
		\begin{equation}\label{normineq2}
			\Big|\lambda_1^*(X_{\beta})-\abs{\lambda_1}\sum_{k\in\beta}u_k^2\Big|\leq \lambda_1^-(X).
		\end{equation} 
		Now we replace $X$ by $G_\alpha$. On the event $B_{\alpha}$ we have
		$$
		\lambda_1^{-}(G_\alpha)\leq \sqrt{m \log \log n} \mbox{ and }\,\,\, y_m<\lambda_1 <y_m+1, 
		$$
		which together with \eqref{normineq2} imply 
		\begin{equation}\label{xbeta}
			y_m \sum_{k\in\beta}u_k^2-\sqrt{m\log \log n}\leq 
			\lambda_1^*(G_{\beta})\leq 
			(y_m+1)\sum_{k\in\beta}u_k^2+\sqrt{m\log \log n}  \end{equation}for all $\beta \subset \alpha$.

		By \eqref{xbeta} and the fact that the first row  $\mathbf{u}$ of the orthogonal group are independent of the eigenvalues  
		$\lambda_i(G_{\alpha}),1\leq i\leq m$, we have 
		\begin{equation}\label{sumxk}
			\begin{split}
				&  \P\left( \sum_{k\in \beta}u_k^2 \leq \frac{y_k-\sqrt{m\log \log n}}{y_m+1}, \forall \beta \subset \alpha, \abs{\beta}=k \right) \\
				\leq & 
				\P(H_{\alpha}|B_{\alpha}) \\
				\leq  &\P\left( \sum_{k\in \beta}u_k^2 \leq \frac{y_k+\sqrt{m\log \log n}}{y_m},\forall \beta \subset \alpha, \abs{\beta}=k \right).
			\end{split}
		\end{equation}
		This together with the facts that $y_k\sim 2\sqrt{k\log n}$  and $y_m\sim 2\sqrt{m\log n}$ as $n\to+\infty$ imply   $$
		\lim_{n\to+\infty}  \P(H_{\alpha}|B_{\alpha})= \P\left(\sum_{k\in \beta} u_k^2\leq \sqrt{k/m}, \forall \beta \in \alpha, \abs{\beta}=k\right),
		$$
		which is exactly the definition of $K_{m}$. This proves \eqref{s4}, and hence \eqref{s1}. Therefore, we complete the proof of Lemma \ref{t_nlimit}.
	\end{proof}
	
	\section{Proof of Theorem \ref{main2}}\label{4}
	Now we prove Theorem \ref{main2}. As before, for   $G=({g_{ij}})_{1\leq i\leq j \leq n}$ sampled from GOE, we denote ${G}_\alpha=({g}_{ij})_{i,j\in \alpha}$ as the principal minor  of size $m\times m$ for $\alpha\subset\{1,..., n\}$ with $|\alpha|=m$ and  $I_m$ is the collection of all such $\alpha$. 
	Let $v_1(G_{\alpha})$ be the eigenvector of the largest eigenvalue of  $G_{\alpha}$ and $v^*\in S^{m-1}$ be the eigenvector of the largest eigenvalue of the principal sub-matrix that attains the maximal eigenvalue $T_{m,n}$, i.e., we have 
	$$
	\alpha^*:=\mathrm{argmax}_{\alpha \in I_m} \lambda_1(G_{\alpha})
	$$
	and
	\begin{equation}
		v^*=v_1(G_{\alpha^*}).
	\end{equation}
	As in \S\ref{lemma44}, we recall the definition of the event 
	\begin{equation}
		H_{\alpha}:=\{\lambda_1^*(G_{\beta})\leq y_k, \, \forall 1\leq k<m,\,\beta \subset \alpha, \abs{\beta}=k \}.
	\end{equation}
	We define the random variable $\hat{\alpha}$ as follows. If the event $\cup_{\alpha \in I_m}H_{\alpha}$ holds, then we set
	$$
	\hat{\alpha}:= \mathrm{argmax}_{\alpha \in I_m, \,\,H_{\alpha}\, \mathrm{ holds}}\lambda_1(G_{\alpha}).
	$$ Otherwise, we set $\hat{\alpha}$ to be $\{1,\ldots, m\}$.  
	We now set
	\begin{equation}
		\hat{v}:=v_1(G_{\hat{\alpha}}).
	\end{equation}
	In other words, $\hat{v}$ is the eigenvector of the largest eigenvalue of  the principal sub-matrix $G_{\alpha}$ that achieves the maximal eigenvalue under the constraint that $H_{\alpha}$ is true. By  \eqref{step1} and \eqref{limit0}, we have   
	$$
	\lim_{n\to+\infty}\P( \cap_{\alpha \in I_m} H_{\alpha})=1. 
	$$
	On the event $\cap_{\alpha \in I_m}H_{\alpha}$ we clearly have $\hat{\alpha}=\alpha^*$  and $\hat{v}=v^*$ since the constraint for $\hat{\alpha}$ doesn't have any effect. 
	Recall the definition of $A_\alpha$ in  \eqref{defAS} where $A_\alpha=H_\alpha\cap \{\lambda_1(G_\alpha)>y_m\}$,  on $\cap_{\alpha \in I_m}H_{\alpha}$,   the two events $\{T_{m,n}\geq y_m \}$ and $\cup_{\alpha \in I_m}A_{\alpha}$ coincide. 
	In other words, for symmetric $Q$, we have 
	\begin{equation}
		\begin{split}
			\left(\cap_{\alpha\in I_m}H_{\alpha} \right) \cap
			\{T_{G,m,n}\geq y_m,v^*\in Q\}  =  \left(\cap_{\alpha\in I_m}H_{\alpha} \right) \cap \left\{\cup_{\alpha\in I_m}A_{\alpha}, \hat{v}\in Q  \right\}.
		\end{split}
	\end{equation}
	It follows that 
	\begin{equation}
		\begin{split}
			&	\abs{	\P(T_{m,n}\geq y_m,v^*\in Q)-\P(\cup_{\alpha\in I_m  } A_{\alpha},\hat{v}\in Q)  }\\	
			=& \abs{  \P\left(  \left( \cap_{\alpha\in I_m}H_{\alpha} \right)^c \cap
				\{T_{m,n}\geq y_m,v^*\in Q \}  \right)
				-\P\left(
				\left(\cap_{\alpha\in I_m}H_{\alpha} \right)^c \cap \left(\cup_{\alpha\in I_m}A_{\alpha}, \hat{v}\in Q  \right)
				\right)} \\
			\leq & 2\P\left(  \left( \cap_{\alpha \in I_m} H_{\alpha} \right)^c  \right),
		\end{split}
	\end{equation}
	which converges to 0 as $n\to+\infty$. Hence, \eqref{value+vec} is equivalent to  the limit
	\begin{equation}\label{value+vec2}
		\P(\cup_{\alpha\in I_m} A_{\alpha}, \hat{v} \in Q ) \to  (1-F_Y(y))\nu (Q).
	\end{equation}
	All of the rest is to prove this convergence. Now we define four quantities 
	\begin{equation}
		\begin{split}
			&	k_1(\alpha)=  \P(A_{\alpha},v_1(G_{\alpha})\in Q; \forall \alpha'\neq \alpha, A_{\alpha'} \mbox{ fails or } \lambda_1(G_{\alpha'})<\lambda_1(G_{\alpha})   ),\\
			& k_2(\alpha)=  \P(A_{\alpha},v_1(G_{\alpha})\in Q; \forall \alpha'\cap \alpha=\emptyset, A_{\alpha'} \mbox{ fails or } \lambda_1(G_{\alpha'})<\lambda_1(G_{\alpha})   ),\\
			&	k_3(\alpha)=  \P(A_{\alpha}; \forall \alpha'\neq \alpha, A_{\alpha'} \mbox{ fails or } \lambda_1(G_{\alpha'})<\lambda_1(G_{\alpha})   ),\\
			& k_4(\alpha)=  \P(A_{\alpha}; \forall \alpha'\cap \alpha=\emptyset, A_{\alpha'} \mbox{ fails or } \lambda_1(G_{\alpha'})<\lambda_1(G_{\alpha})   ).
		\end{split}
	\end{equation}
	We note that
	\begin{equation}\label{555}
		\begin{split}
			&\P\left(   \cup_{\alpha\in I_m} A_{\alpha},\hat{v}\in Q \right)\\
			=&\sum_{\alpha\in I_m}\P(A_{\alpha}, v_1(G_{\alpha})\in Q ; \forall \alpha'\neq \alpha, A_{\alpha'} \mbox{ fails or } \lambda_1(G_{\alpha'})<\lambda_1(G_{\alpha})   )\\
			=&\sum_{\alpha \in I_m}k_1(\alpha).
		\end{split}
	\end{equation}
	In fact $k_1(\alpha),\ldots,k_4(\alpha)$ don't depend on the specific choice of $\alpha$. 
	Recall the definition of $b_{n,2}$ and \eqref{bn2ctl}, we have 
	\begin{equation}\label{2bn2}
		\lim_{n\to+\infty}\sum_{\alpha,\alpha' \in I_m;\alpha\cap \alpha' \neq \emptyset}  \P(A_{\alpha}\cap A_{\alpha'}) =\lim_{n\to+\infty}b_{n,2}=0.
	\end{equation}
	By the definition of $k_1(\alpha)$ and $k_2(\alpha )$ together with the  union bound we have 
	\begin{equation}\label{k1k2}
		\sum_{\alpha \in I_m}	\abs{k_1(\alpha)-k_2(\alpha)}\leq \sum_{\alpha\in I_m} \P\left(A_{\alpha} \cap \left( \cup_{\alpha' \cap \alpha \neq \emptyset} A_{\alpha'}\right)  \right) \leq 
		\sum_{\alpha,\alpha' \in I_m;\alpha\cap \alpha' \neq \emptyset}  \P(A_{\alpha}\cap A_{\alpha'}).
	\end{equation}
	Combining \eqref{555}, \eqref{2bn2} and \eqref{k1k2} we see that 
	\begin{equation}\label{2k2}
		\lim_{n\to+\infty} \abs{\P(\cup_{\alpha \in I_m} A_{\alpha},\hat{v}\in Q ) -\sum_{\alpha\in I_m}k_2(\alpha)}=0. 
	\end{equation}
	We can similarly show that
	\begin{equation}\label{k2app}
		\begin{split}
			\lim_{n\to+\infty}\abs{\P(\cup_{\alpha\in I_m} A_{\alpha}) -\sum_{\alpha\in I_m} k_4(\alpha)}
			\leq \lim_{n\to+\infty}\sum_{\alpha\in I_m}\abs{k_3(\alpha)-k_4(\alpha)}	=0.	
		\end{split}
	\end{equation}
	Recall that we have already shown 
	\begin{equation}\label{k4app}
		\P(\cup_{\alpha\in I_m} A_{\alpha}) \to 1-F_Y(y)
	\end{equation}
	in the proof of Theorem \ref{main} as $n\to+\infty$. Hence, combining \eqref{k2app} and \eqref{k4app}  we have
	\begin{equation}\label{k4limit}
		\lim_{n\to+\infty} \sum_{\alpha\in I_m} k_4(\alpha)=1-F_Y(y). 
	\end{equation}
	The advantage of introducing $k_2(\alpha)$ is that, we have the conditional probability 
	\begin{equation}
		\begin{split}
			&	\P(v_1(G_{\alpha})\in Q|A_{\alpha};\forall \alpha'\cap \alpha=\emptyset, A_{\alpha'} \mbox{ fails or }\lambda_1(G_{\alpha'})<\lambda_1(G_{\alpha}) )\\
			=& \P(v_1(G_{\alpha})\in Q|A_{\alpha}), 
		\end{split}
	\end{equation}
	since $G_{\alpha}$ is independent of $\{G_{\alpha'}:\alpha'\cap \alpha =\emptyset\}$.
	Consequently, we have
	\begin{equation}
		\begin{split}
			k_2(\alpha)= 	&	
			\P(A_{\alpha};\forall \alpha'\cap \alpha=\emptyset, A_{\alpha'} \mbox{ fails or }\lambda_1(G_{\alpha'})<\lambda_1(G_{\alpha}) ) \\
			& \times 	\P(v_1(G_{\alpha})\in Q|A_{\alpha};\forall \alpha'\cap \alpha=\emptyset, A_{\alpha'} \mbox{ fails or }\lambda_1(G_{\alpha'})<\lambda_1(G_{\alpha}) )\\
			=&\P(v_1(G_{\alpha})\in Q|A_{\alpha}) k_4(\alpha) .
		\end{split}
	\end{equation}
	Clearly the term  $\P(v_1(G_{\alpha})\in Q|A_{\alpha})$ is the same for all $\alpha \in I_m$. 
	Let $\alpha_1=\{1,\ldots,m\}$. Then we have
	\begin{equation}\label{k2sum}
		\sum_{\alpha \in I_{m}}k_2(\alpha)=\sum_{\alpha \in I_m}  \P(v_1(G_{\alpha})\in Q|A_{\alpha}) k_4(\alpha)=\P(v_1(G_{\alpha_1})\in Q |A_{\alpha_1})\sum_{\alpha\in I_m} k_4(\alpha).
	\end{equation}
	
	Let $\mathbf u$ be sampled from the uniform distribution on the unit sphere $S^{m-1}$. 
	Inspecting the proof of \eqref{s1} in \S\ref{lemma44}, especially \eqref{s2} and \eqref{sumxk},  we have  
	\begin{equation}\label{v1qq}
		\lim_{n\to+\infty}\P(H_{\alpha},v_1(G_{\alpha}) \in Q|\lambda_1(G_{\alpha})>y_m )=\P({\mathbf u}\in \mathcal S_m \cap Q), 
	\end{equation}
	where  $\mathcal S_m$ has been defined in \eqref{hm}.
	
	By the fact $A_{\alpha}=H_{\alpha}\cap \{\lambda_1(G_{\alpha})>y_m \},$ we have
	\begin{equation}\label{v1limit}\begin{split}		
			&\lim_{n\to+\infty}	\P(v_1(G_{\alpha_1})\in Q|A_{\alpha_1})\\=&
			\lim_{n\to+\infty}   \frac{\P(H_{\alpha},v_1(G_{\alpha}) \in Q|\lambda_1(G_{\alpha})>y_m )}{\P(H_{\alpha}|\lambda_1(G_{\alpha})>y_m )}\\=&\frac{\P({\mathbf u}\in \mathcal S_m \cap Q)}{\P({\mathbf u}\in \mathcal S_m)}\\=&
			\nu(Q), \end{split}
	\end{equation}
	where  $\nu$ is the uniform distribution on the set $\mathcal S_m$. 
	
	Combining \eqref{2k2}, \eqref{k4limit}, \eqref{k2sum} and \eqref{v1limit}  we get 
	\begin{equation}
		\begin{split}
			&\lim_{n\to+\infty} \P\left(\cup_{\alpha\in I_m} A_{\alpha},\hat{v}\in Q \right)\\
			=& \lim_{n\to+\infty} \sum_{\alpha\in I_m}k_2(\alpha)\\
			=& \lim_{n\to+\infty}  	\P(v_1(G_{\alpha_1})\in Q|A_{\alpha_1}) \sum_{\alpha\in I_m}k_4(\alpha)\\
			=& \nu(Q) (1-F_Y(y)).
		\end{split}
	\end{equation}
	This proves \eqref{value+vec2},  and thus    the proof of Theorem \ref{main2}.


\begin{thebibliography}{99}
		\bibitem{AGZ}G. W. Anderson,  A. Guionnet and O. Zeitouni, \emph{An introduction to random matrices}. Cambridge Studies in Advanced Mathematics, 118. Cambridge University Press, Cambridge, 2010.
		
		\bibitem{AGG}Arratia R., Goldstein L., Gordon L. Two moments suffice for Poisson approximations:
		The Chen-Stein method. Ann. Probab. 17(1989):9-25.
		
		\bibitem{B}Baraniuk, R., Davenport, M., DeVore, R., and Wakin, M. (2008). A simple proof of the
		restricted isometry property for random matrices. Constructive Approximation, 28(3):253-263.
		
		
		
		
		\bibitem{CJL}Cai T.-T., Jiang T.,  Li X. Asymptotic Analysis for Extreme Eigenvalues of
		Principal Minors of Random Matrices, Ann. Appl. Probab. 31 (6) 2953 - 2990, December 2021.
		
		
		
		\bibitem{GL} Gamarnik, D., and Li, Q. (2018). Finding a large submatrix  of a Gaussian Random matrix. The Annals of Statistics, 46(6A), 2511-2561.
		
		
		
		\bibitem{LLr} Leadbetter, M.R., Lindgren, G., Rootzen, H.: \emph{Extremes and related properties of random
			sequences and processes}. Springer Series in Statistics, Springer 
		(1983).
		
		\bibitem{LR}	Leadbetter, M. R. and Rootzen, H. (1988). Extremal theory for stochastic processes. Ann.  Probab., 431-478.
		
		
		
		
		
		
		
		
		\bibitem{TC2}Tao, Terence; Cand\`{e}s, Emmanuel J. (2006), Decoding by linear programming, IEEE transactions on information theory 51 (12), 4203-4215. 
		
		
\end{thebibliography}
	 \end{document}